\documentclass[onefignum,onetabnum]{siamart190516}

\usepackage{tikz}
\usetikzlibrary{scopes}
\usetikzlibrary{arrows}
\usetikzlibrary{arrows.meta}
\usetikzlibrary{decorations.pathreplacing,decorations.markings}
\usepackage{lipsum}
\usepackage{amsfonts}
\usepackage{amsmath}
\usepackage{graphicx}
\usepackage{subfigure}
\usepackage{epstopdf}
\usepackage{algorithmic}
\ifpdf
\DeclareGraphicsExtensions{.eps,.pdf,.png,.jpg}
\else
\DeclareGraphicsExtensions{.eps}
\fi

\newsiamremark{remark}{Remark}
\newsiamremark{hypothesis}{Hypothesis}
\crefname{hypothesis}{Hypothesis}{Hypotheses}
\newsiamthm{claim}{Claim}

\headers{Numerical procedure for hybrid systems}{R. Pytlak and D. Suski}

\title{Numerical procedure for optimal control of hybrid systems with sliding modes, Part II}

\author{Rados\L aw Pytlak\thanks{Faculty of Mathematics and Information Science,
		Warsaw University of Technology, 00-665 Warsaw, Poland 
		(\email{r.pytlak@mini.pw.edu.pl}),}
	\and Damian Suski\thanks{Institute of Automatic Control and Robotics,
		Warsaw University of Technology, 02-525 Warsaw, Poland 
		(\email{d.suski@mchtr.pw.edu.pl}).}
}

\usepackage{amsopn}

\begin{document}

\maketitle

\begin{abstract}
This paper concerns the numerical procedure for solving hybrid optimal control problems with sliding modes. 
A sliding mode is coped with differential-algebraic equations (DAEs) and that guarantees accurate tracking of the sliding motion surface. 
In the second part of the paper we demonstrate the correspondence between the discrete adjoint equations and the discretized version of the continuous adjoint equations in the case of system equations described by DAEs. We show that the discrete adjoint state trajectories converge to their continuous counterparts.  
Next, we describe the application of the proposed procedure to three optimal control problems. The first problem concerns optimal control of a simple mechanical system with dry friction. The second problem is related to the planning of a haemodialysis process. The third problem concerns the optimal steering of a racing car. 
\end{abstract}



\section{Introduction}
\label{secIntro}

The second part of the paper continues the description of the numerical procedure for optimal control problems for hybrid system exhibiting sliding motion behavior. The procedure is the implementation of the algorithm for these problems introduced in \cite{ps19}. In the first part of the paper we concentrated on describing the efficient way of evaluating reduced gradients of the optimization problem. We paid the attention to the case the hybrid system does not enter the sliding mode and its dynamics is described by ODEs. In the second part we supplement our considerations by analyzing the sliding motion case. 

Let us consider the hybrid system evolving according to the following scheme. For $t\in \left[t_0,t_t\right]$ the system evolves according to ODEs
\begin{equation}
x' = f^1(x,u),
\label{eqAdjEqsFirstCaseOde}
\end{equation}
at a switching time $ t_t $ system state reaches the switching surface 
\begin{equation}
g(x(t_t)) = 0
\label{eqAdjEqsFirstCaseTrans}
\end{equation} 
and enters the sliding mode. For  $t\in \left[t_t,t_f\right]$ the state trajectory is the solution of DAEs 
\begin{eqnarray}
x' &=& f_F(x,u) + g_x^T(x)z = f^2(x,z,u) \label{eqAdjEqsFirstCaseDaeDiff} \\
0 &=& g(x) \label{eqAdjEqsFirstCaseDaeAlg}
\end{eqnarray}

We now aim at finding the reduced gradient formula for the endpoint functional (\cite{ps19a})
\begin{equation}
\bar{F}_0(u) = \phi(x^u(t_f))
\end{equation}
with respect to controls.
This formula can be found with the help of adjoint equations. We now recall the appropriate formulation of adjoint equations in the considered case. The consistent terminal values of the adjoint variables $\lambda_f(t_f)$ and $\lambda_g(t_f)$ can be found by solving the following set of equations with respect to the variables $ \lambda_f(t_f), \lambda_g(t_f), \nu_1 $ (\cite{ps19})
\begin{eqnarray}
\phi_x^T(x(t_f)) + \lambda_f(t_f) & = & \nu_1 g_x^T(x(t_f)) \label{terminal} \\
0 & = & g_x(x(t_f))\lambda_f(t_f)  \\
0 & = & \left(g_x(x(t_f))\right)'\lambda_f(t_f) - \label{slidingDAEIndex2HiddenAdjEndpoint} \\
& & g_x(x(t_f)) \left( f_F \right)_x^T(x(t_f),u(t_f))\lambda_f(t_f) -\nonumber \\
& & g_x(x(t_f)) \left(g_x^T(x(t_f))z(t_f)\right)_x^T\lambda_f(t_f) +  \nonumber \\
& & g_x(x(t_f)) g_x^T(x(t_f)) \lambda_g (t_f) \nonumber 
\end{eqnarray}
where $ \nu_1 $ is some real number.

Having $\lambda_f(t_f)$ and $\lambda_g(t_f)$ we then solve DAEs (the meaning of $t_t^-$ and $t_t^+$ is explained in \cite{ps19a})
\begin{eqnarray}
(\lambda_f^T)'(t) &= &-\lambda_f^T(t) (f_F)_x(x(t),u(t)) \label{eqAdjEqsFirstCaseAdjDaeDiff}  \\
&\ & - \lambda_f^T(t) (g_x^T(x(t))z(t))_x + \lambda_g^T(t) g_x(x(t)) \nonumber \\
0 &= & \lambda_f^T(t) g_x^T(x(t)),\ t\in [t_t^+,t_f]. \label{eqAdjEqsFirstCaseAdjDaeAlg}
\end{eqnarray}
backwards in time. At the transition time $ t_t $ the adjoint variable $ \lambda_f $ undergoes a jump. To calculate the value of $ \lambda_f(t_t^-) $ the following system of equations have to be solved for the variables $ \lambda_f(t_t^-),\ \pi $ 
\begin{eqnarray}
\lambda_f(t_t^-) & = & \lambda_f(t_t^+) - \pi g_x^T(x(t_t))  \label{odeAdjTransHamilConta} \\
\lambda_f^T(t_t^-) f^1(x(t_t^-),u(t_t^-)) & = & \lambda_f^T(t_t^+) f_F(x(t_t^+),u(t_t^+)) +  \label{odeAdjTransHamilCont}  \\
& & \lambda_f^T(t_t^+)g_x^T(x(t_t^+))z(t_t^+) - \lambda_g^T(t_t^+)g(x(t_t^+)). \nonumber
\end{eqnarray}
Eventually we solve adjoint ODEs  
\begin{equation}
(\lambda_f^T)'(t) = - \lambda_f^T(t) (f^1)_x(x(t),u(t)),\ t\in[t_0,t_t^-).  
\label{eqAdjEqsFirstCaseAdjOde}
\end{equation}

Now the directional derivative of the functional $ \bar{F}_0(u) $ with respect to the control variation $ d $ can be calculated from
\begin{eqnarray}
\langle \nabla \bar{F}_0(u),d\rangle  &= &-\int_{t_0}^{t_t^-} \lambda_f^T(t) (f^1)_u(x(t),u(t))d(t) dt \label{eqAdjGradient}\\
& &   -\int_{t_t^+}^{t_f} \lambda_f^T(t) (f_F)_u(x(t),u(t))d(t) dt. \nonumber  
\end{eqnarray} 

In the next section (Section \ref{secDiscrete}) we will analyze the discrete time versions of equations (\ref{eqAdjEqsFirstCaseAdjDaeDiff})--(\ref{eqAdjEqsFirstCaseAdjDaeAlg}) and the approximation to the second term in (\ref{eqAdjGradient}). The present analysis complements the analysis given in the first part of the paper \cite{ps19a} which concerned equations (\ref{eqAdjEqsFirstCaseAdjOde}) and the first term in (\ref{eqAdjGradient}). In addition in Section \ref{secJump} we will examine how the discretization of adjoint equations influences the jump conditions (\ref{odeAdjTransHamilCont}). Section \ref{secResults} presents the efficiency of the proposed numerical procedure in solving three optimal control problems.

\section{Numerical calculation of reduced gradients} \label{secDiscrete}

Let us consider the control system described by higher index differential-algebraic equations, we focus on index 2 DAEs in Hessenberg form (\cite{hlr89})
\begin{eqnarray}
x'(t) & = & f(x(t),z(t),u(t)), \label{discrEq:daeDiff}\\
0 & = & g(x(t)), \label{discrEq:daeAlg}
\end{eqnarray} 
where $ x\in\mathbb{R}^{n_x} $ is a differential state and $ z\in\mathbb{R}^{n_z} $ is an algebraic state. The differentiation index of a system (\ref{discrEq:daeDiff})-(\ref{discrEq:daeAlg}) is 2 provided that the matrix $ g_x(x)f_z(x,z,u) $ is nonsingular (\cite{hlr89}). To integrate DAEs, the following Runge-Kutta scheme can be used (\cite{hlr89})
\begin{eqnarray}
x_i'(k+1) &=& f \left( x_i(k+1),\ z_i(k+1), u(k) \right),  \label{discrEq:rkSchemeIntDiff001} \\
0 &=& g \left( x_i(k+1)  \right), \label{discrEq:rkSchemeIntAlg001} \\
x_i(k+1) &=& x(k)+h(k)\sum_{j=1}^{s}a_{ij}x_j'(k+1), \label{discrEq:rkSchemeIntDiff002} \\
z_i(k+1) &=& z(k)+h(k)\sum_{j=1}^{s}a_{ij}z_j'(k+1) \label{discrEq:rkSchemeIntAlg002} 
\end{eqnarray}
for $ i = 1,...,s $ and 
\begin{eqnarray}
x(k+1) &=& x(k) +h(k)\sum_{i=1}^{s}b_i x_i'(k+1), \label{discrEq:rkSchemeStepDiff001} \\
z(k+1) &=& z(k) +h(k)\sum_{i=1}^{s}b_i z_i'(k+1). \label{discrEq:rkSchemeStepAlg001} 
\end{eqnarray}
It is possible to get rid of variables $ x_i'(k+1) $ and $ z_i'(k+1) $ from (\ref{discrEq:rkSchemeIntDiff001})-(\ref{discrEq:rkSchemeStepAlg001}), provided that the Runge-Kutta matrix $ A = (a_{ij}) $ is invertible. In such case, from (\ref{discrEq:rkSchemeIntAlg002}) we have
\begin{equation}
z_i'(k+1) = \frac{1}{h(k)} \sum_{j=1}^{s} a_{ij}^- (z_j(k+1) - z(k)),\ i = 1,\dots,s, 
\end{equation}  
where $ a_{ij}^- $ are the coefficients of the matrix $ A^{-1} $. Now the system (\ref{discrEq:rkSchemeIntDiff001})-(\ref{discrEq:rkSchemeStepAlg001}) can be rewritten as
\begin{eqnarray}
x_i(k+1) &=& x(k)+h(k)\sum_{j=1}^{s}a_{ij}f \left( x_j(k+1),\ z_j(k+1), u(k) \right), \label{discrEq:rkSchemeIntDiff003} \\
0 &=& g \left( x_i(k+1)  \right), \label{discrEq:rkSchemeIntAlg003} 
\end{eqnarray}
for $ i = 1,...,s $ and 
\begin{eqnarray}
x(k+1) &=& x(k) +h(k)\sum_{i=1}^{s}b_i f \left( x_i(k+1),\ z_i(k+1), u(k) \right), \label{discrEq:rkSchemeStepDiff002} \\
z(k+1) &=& z(k) +\sum_{i=1}^{s} \sum_{j=1}^{s} b_i  a_{ij}^- (z_j(k+1) - z(k)). \label{discrEq:rkSchemeStepAlg002} \\
&=& z(k) +\sum_{i=1}^{s}  b_i^-  (z_i(k+1) - z(k)), \nonumber
\end{eqnarray} 
where
\begin{equation}
b_i^- = \sum_{i=1}^{s} b_j a^-_{ji},\ i = 1,\dots,s.
\end{equation} 
At each step of the Runge--Kutta scheme the nonlinear system (\ref{discrEq:rkSchemeIntDiff003})-(\ref{discrEq:rkSchemeIntAlg003}) is first solved for variables $ x_i(k+1), z_i(k+1)\ i = 1,\dots,s $ and then (\ref{discrEq:rkSchemeStepDiff002})-(\ref{discrEq:rkSchemeStepAlg002}) are used to calculate $ x(k+1) $ and $ z(k+1) $. Let us notice that $ z(k) $ is required to calculate $ z(k+1) $, but not for solving (\ref{discrEq:rkSchemeIntDiff003})-(\ref{discrEq:rkSchemeIntAlg003}).

The adjoint equations for the control system (\ref{discrEq:daeDiff})-(\ref{discrEq:daeAlg}) are (\cite{ps19}):
\begin{eqnarray}
\lambda_f '(t) &=& -f_x^T(x(t), z(t), u(t)) \lambda_f(t) - g_x^T(x(t))\lambda_g(t), \label{discrEq:adjDAEDiff} \\
0 &=& - f_z^T(x(t), z(t), u(t)) \lambda_f(t). \label{discrEq:adjDAEAlg}
\end{eqnarray}
For the adjoint DAEs (\ref{discrEq:adjDAEDiff})-(\ref{discrEq:adjDAEAlg}) to be properly defined, the term $ f_z^T(x, z, u)\lambda_f $ should be differentiable. The algebraic state $ z(t) $ and the control $ u(t) $ are not continuous, so to keep the differentiability we assume that $  f_z^T(x(t), z(t), u(t)) $ does not depend neither on $ z(t) $ nor on $ u(t) $. Let us notice that the adjoint equations formulated for the sliding motion satisfy that condition. From now we will consider the adjoint equations of the form
\begin{eqnarray}
\lambda_f '(t) &=& -f_x^T(x(t), z(t), u(t)) \lambda_f(t) - g_x^T(x(t))\lambda_g(t), \label{discrEq:adjDAEDiff001} \\
0 &=& - f_z^T(x(t)) \lambda_f(t). \label{discrEq:adjDAEAlg001}
\end{eqnarray}
The adjoint equations are index 2 DAEs (\cite{ps19}), so the Runge-Kutta scheme can be used to numerically integrate them (\cite{hlr89}). The Runge-Kutta scheme for adjoint equations (\ref{discrEq:adjDAEDiff001})-(\ref{discrEq:adjDAEAlg001}) is
\begin{eqnarray}
&&\lambda_{fi}(k)= \nonumber \\ && \lambda_f(k+1) - \bar{h}(k+1) \sum_{j=1}^{\bar{s}} \bar{a}_{ij} \left[  -f_x^T 
\left( 
x\left( \bar{t}_j(k) \right) , z\left( \bar{t}_j(k) \right), \bar{u}(k+1) 
\right) \lambda_{fj}(k) \right.
\nonumber\\
&&\qquad\qquad\qquad\qquad\qquad\qquad\ \
\left.  -g_x^T 
\left( 
x\left( \bar{t}_j(k) \right) \right) \lambda_{gj}(k) \right] =  \label{discrEq:rkSchemeIntAdjDiff001} \\
&& \lambda_f(k+1) +  \bar{h}(k+1) \sum_{j=1}^{\bar{s}} \bar{a}_{ij} \left[  f_x^T 
\left( 
x\left( \bar{t}_j(k) \right) , z\left( \bar{t}_j(k) \right), \bar{u}(k+1) 
\right) \lambda_{fj}(k) \right.
\nonumber\\
&&\qquad\qquad\qquad\qquad\qquad\qquad\ \
\left.  +g_x^T 
\left( 
x\left( \bar{t}_j(k) \right) \right) \lambda_{gj}(k) \right], \nonumber \\
&& 0 = -f^T_z\left( x\left( \bar{t}_i(k) \right) \right) \lambda_{fi}(k)\label{discrEq:rkSchemeIntAdjAlg001}
\end{eqnarray}
for $ i = 1,\dots,\bar{s} $ and
\begin{eqnarray}
&&\lambda_{f}(k) = \nonumber \\ 
&&\lambda_f(k+1) - \bar{h}(k+1) \sum_{i=1}^{\bar{s}} \bar{b}_{i} \left[  -f_x^T 
\left( 
x\left( \bar{t}_i(k) \right) , z\left( \bar{t}_i(k) \right), \bar{u}(k+1) 
\right) \lambda_{fi}(k) \right.
\nonumber\\
&&\qquad\qquad\qquad\qquad\qquad\qquad\ \
\left.  -g_x^T 
\left( 
x\left( \bar{t}_i(k) \right) \right) \lambda_{gi}(k) \right]=  \label{discrEq:rkSchemeStepAdjDiff001} \\
&&\lambda_f(k+1) + \bar{h}(k+1) \sum_{i=1}^{\bar{s}} \bar{b}_{i} \left[  f_x^T 
\left( 
x\left( \bar{t}_i(k) \right) , z\left( \bar{t}_i(k) \right), \bar{u}(k+1) 
\right) \lambda_{fi}(k) \right.
\nonumber\\
&&\qquad\qquad\qquad\qquad\qquad\qquad\ \
\left.  +g_x^T 
\left( 
x\left( \bar{t}_i(k) \right) \right) \lambda_{gi}(k) \right], \nonumber \\
&&\lambda_{g}(k) = \nonumber \\
&&\lambda_{g}(k+1) +\sum_{i=1}^{\bar{s}} \sum_{j=1}^{\bar{s}} \bar{b}_i  \bar{a}_{ij}^- (\lambda_{gj}(k) - \lambda_{g}(k+1) ) = \label{discrEq:rkSchemeStepAdjAlg001} \\
&&\lambda_{g}(k+1) +\sum_{i=1}^{\bar{s}} \bar{b}_i^- (\lambda_{gi}(k) - \lambda_{g}(k+1) ), \nonumber 
\end{eqnarray}
where for the sake of the shorter notation we denote
\begin{equation}
\bar{t}_i(k) = t(k+1) - \bar{c}_i \bar{h}(k+1),\ i = 1,\dots,s.
\end{equation}
At each step of the Runge-Kutta scheme, the system of equations (\ref{discrEq:rkSchemeIntAdjDiff001})-(\ref{discrEq:rkSchemeIntAdjAlg001}) is first solved for variables $ \lambda_{fi}(k),\lambda_{gi}(k),\ i = 1,\dots,\bar{s} $ and then (\ref{discrEq:rkSchemeStepAdjDiff001})-(\ref{discrEq:rkSchemeStepAdjAlg001}) are used to calculate $ \lambda_f(k) $ and $ \lambda_g(k) $.

In the case of the sliding motion the differential of the endpoint functional $\phi(t_f)$ is of the form 
\begin{equation}
d \phi(x(t_f)) = \int_{t_0}^{t_f} d(t)^T f_u^T(x(t),u(t)) \lambda_f(t) dt.
\label{discrEq:funcDiffContSlid}
\end{equation}
Two important remarks need to be emphasized. First we assume that $ f_u(x,u) $ does not depend on $ z $. That is another limitation lied on the form of the function $ f(x,z,u) $, but it is satisfied by the sliding motion equations. Second, the integral depends only on the adjoint variable $ \lambda_{f} $, but not on $ \lambda_{g} $. Under that conditions the reduced gradients for DAEs are calculated the same way as it was for ODEs (see \cite{ps19a})
\begin{equation}
\frac{d \phi(x(t_f))}{ du_n } \simeq \sum_{k=k_{n-1}}^{k_n-1} \tilde{h}(k) \sum_{i=1}^{\tilde{s}} \tilde{b}_i f_u^T \left(x \left(\tilde{t}(k)+\tilde{c}_i \tilde{h}(k)\right),u_n \right) \lambda_f \left( \tilde{t}(k)+\tilde{c}_i \tilde{h}(k) \right),
\label{discrEq:redGradQuadContSlid}
\end{equation}
where for the discrete steps $ k = k_{n-1},\hdots, k_n-1 $ the control is $ u_n $.

Similarly to the ODEs case presented in \cite{ps19a}, we utilize discrete adjoint equations and discrete reduced gradients to avoid the necessity of calculation of states and adjoint variables at arbitrary time moments. The Runge-Kutta scheme (\ref{discrEq:rkSchemeIntDiff003})-(\ref{discrEq:rkSchemeStepDiff002}) can be rewritten to a vector form (discrete step argument omitted)
\begin{equation}
\left( 
\begin{array}{c}
x_1 - x - h\sum_{j=1}^{s}a_{1j}f \left( x_j, z_j, u \right) \\
- g \left( x_1  \right)\\
\vdots \\
x_s - x - h\sum_{j=1}^{s}a_{sj}f \left( x_j, z_j, u \right) \\
- g \left( x_s \right) \\
x^+ - x - h\sum_{i=1}^{s}b_i f \left( x_i, z_i, u \right)
\end{array}
\right) = 
\left( 
\begin{array}{c}
0 \\ 
0 \\ 
\vdots \\ 
0 \\ 
0 \\ 
0 
\end{array}
\right).
\label{discrEq:rkSchemeVecDAE}
\end{equation}
If we now define the augmented state vector $ X(k) $ as
\begin{equation}
X(k) = \left( x_1(k)^T,z_1(k)^T,\ldots,x_s(k)^T,z_s(k)^T, x(k)^T \right)^T ,
\end{equation}
then (\ref{discrEq:rkSchemeVecDAE}) can be presented in a form of the implicit discrete time state equation \cite{py11} (see also \cite{pz14})
\begin{equation}
F\left(X(k+1),X(k),u(k)\right) = 0.
\label{discrEq:discrStateEq}
\end{equation}
The partial derivatives matrices of the discrete state equation are
\begin{equation}
F_{X^+}(k) = \left(
\begin{array}{ccccccc}
I - ha_{11}f_{x1} & -ha_{11}f_{z1} & \ldots & - ha_{1s}f_{xs} & -ha_{1s}f_{zs} & 0  \\
-g_{x1} & 0 & \ldots & 0 & 0 & 0   \\
\vdots & \vdots & & \vdots & \vdots & \vdots  \\
- ha_{s1}f_{x1} & -ha_{s1}f_{z1} & \ldots & I - ha_{ss}f_{xs} & -ha_{ss}f_{zs} & 0   \\
0 & 0 & \ldots & -g_{xs} & 0 & 0  \\
-h b_1 f_{x1} & -h b_1 f_{z1} & \ldots & -h b_s f_{xs} & -h b_s f_{zs} & I 
\end{array}
\right )\label{discrEq:FXPlusDAE}
\end{equation}
and
\begin{equation}
F_{X}(k) = \left(
\begin{array}{ccccccc}
0 & 0 & \ldots & 0 & 0 & -I    \\
0 & 0 & \ldots & 0 & 0 & 0     \\
\vdots & \vdots & & \vdots & \vdots & \vdots  \\
0 & 0 & \ldots & 0 & 0 & -I    \\
0 & 0 & \ldots & 0 & 0 & 0     \\
0 & 0 & \ldots & 0 & 0 & -I 
\end{array}
\right ),\label{discrEq:FXDAE}
\end{equation}
where for the sake of a shorter notation we denote 
\begin{equation}
f_{xi} = f_{x}(x_i,z_i,u), f_{zi} = f_{z}(x_i), g_{xi} = g_{x}(x_i).
\end{equation}
Let us denote 
\begin{eqnarray}
\Lambda(k) & = & \left( l_{f1}(k)^T,l_{g1}(k)^T, \ldots, l_{fs}(k)^T,l_{gs}(k)^T, \lambda_f(k)^T \right)^T,\nonumber \\
\Lambda(k+1) & = & \left( l_{f1}^+(k)^T,l_{g1}^+(k)^T, \ldots, l_{fs}^+(k)^T,l_{gs}^+(k)^T, \lambda_f^+(k)^T \right)^T,\nonumber \\
R(k) & = & \left( r_{f1}(k)^T,r_{g1}(k)^T, \ldots, r_{fs}(k)^T,r_{gs}(k)^T, r_f(k)^T \right)^T.\nonumber
\end{eqnarray}
Now the discrete adjoint equations (\cite{py11})
\begin{eqnarray}
F_{X^+}^{T} (k) R(k) &=& \Lambda(k+1) \label{discrEq:discrAdjEqRCalc}, \\
\Lambda(k) &=& -F_X^T(k) R(k) \label{discrEq:discrAdjEqLamCalc},
\end{eqnarray}  
at a discrete time step $ k $  take the form
\begin{equation}
\left(
\begin{array}{cccccc}
I - ha_{11}f^T_{x1} & -g^T_{x1} & \ldots & - ha_{s1}f^T_{x1} & 0 & -h b_1 f_{x1}^T  \\
-ha_{11}f^T_{z1} & 0 & \ldots & -ha_{s1}f^T_{z1} & 0 & -h b_1 f_{z1}^T   \\
\vdots & \vdots & & \vdots & \vdots & \vdots  \\
-ha_{1s}f^T_{xs} & 0 & \ldots & I - ha_{ss}f^T_{xs} & -g^T_{xs} & -h b_s f_{xs}^T  \\
-ha_{1s}f^T_{zs} & 0 & \ldots & -ha_{ss} f^T_{zs}& 0 & -h b_s f_{zs}^T    \\
0 & 0 & \ldots & 0 & 0 & I \\	
\end{array}
\right ) 
\left( 
\begin{array}{c}
r_{f1} \\ r_{g1} \\ \vdots \\ 
r_{fs} \\ r_{gs} \\ r_f 
\end{array}
\right) = 
\left(
\begin{array}{c}
l^+_{f1} \\ l^+_{g1} \\ \vdots \\ 
l^+_{fs} \\ l^+_{gs} \\ \lambda^+_f 
\end{array}
\right) \label{discrEq:discrAdjEqDAERCalcVec}
\end{equation}

\begin{equation}
\left(
\begin{array}{c}
l_{f1} \\ l_{g1} \\ \vdots \\ 
l_{fs} \\ l_{gs} \\ \lambda_f 
\end{array}
\right) = -
\left(
\begin{array}{cccccc}
0 & 0 & \ldots & 0 & 0 & 0  \\
0 & 0 & \ldots & 0 & 0 & 0     \\
\vdots & \vdots & & \vdots & \vdots & \vdots  \\
0 & 0 & \ldots & 0 & 0 & 0   \\
0 & 0 & \ldots & 0 & 0 & 0     \\
-I & 0 & \ldots & -I & 0  & -I 
\end{array}
\right)
\left( 
\begin{array}{c}
r_{f1} \\ r_{g1} \\ \vdots \\ 
r_{fs} \\ r_{gs} \\ r_f 
\end{array}
\right) \label{discrEq:discrAdjEqDAELamCalcVec}
\end{equation}
If we rewrite (\ref{discrEq:discrAdjEqDAERCalcVec}) as a system of equations we obtain
\begin{eqnarray}
r_{fi} &=& \sum_{j=1}^{s} ha_{ji} f_{xi}^T r_{fj} + g_{xi}^T r_{gi}+  hb_if_{xi}^T r_f + l^+_{fi} \label{discrEq:discrAdjEqRIntDiffCalcSysEq001} \\
0 &=& \sum_{j=1}^{s} ha_{ji}  f_{zi}^T r_{fj} + hb_if_{zi}^T r_f + l_{gi}^+  \label{discrEq:discrAdjEqRIntAlgCalcSysEq001}
\end{eqnarray}
for $ i = 1,\ldots,s $ and
\begin{equation} 
r_f = \lambda_f^+ \label{discrEq:discrAdjEqREndDiffCalcSysEq001} 
\end{equation}
Using (\ref{discrEq:discrAdjEqREndDiffCalcSysEq001}), (\ref{discrEq:discrAdjEqRIntDiffCalcSysEq001})-(\ref{discrEq:discrAdjEqRIntAlgCalcSysEq001}) can be written as
\begin{eqnarray}
r_{fi} &=& \sum_{j=1}^{s} ha_{ji} f_{xi}^T r_{fj} + g_{xi}^T r_{gi}+  hb_if_{xi}^T \lambda_f^+ + l^+_{fi} \label{discrEq:discrAdjEqRIntDiffCalcSysEq002} \\
0 &=& \sum_{j=1}^{s} ha_{ji}  f_{zi}^T r_{fj} + hb_if_{zi}^T \lambda_f^+ + l_{gi}^+  \label{discrEq:discrAdjEqRIntAlgCalcSysEq002}
\end{eqnarray}

If we rewrite (\ref{discrEq:discrAdjEqDAELamCalcVec}) as a system of equations we get
\begin{eqnarray}
l_{fi} &=& 0 \label{discrEq:discrAdjEqAdjIntDiffCalcSysEq002}\\
l_{gi} &=& 0 \label{discrEq:discrAdjEqAdjIntAlgCalcSysEq002}
\end{eqnarray}
for $ i = 1,\ldots,s $ and
\begin{eqnarray}
\lambda_f &=& \sum_{i=1}^{s} r_{fi} + r_f. \label{discrEq:discrAdjEqAdjEndDiffCalcSysEq001}  
\end{eqnarray}
Using (\ref{discrEq:discrAdjEqREndDiffCalcSysEq001}), (\ref{discrEq:discrAdjEqAdjEndDiffCalcSysEq001}) gives 
\begin{eqnarray}
\lambda_f &=& \sum_{i=1}^{s} r_{fi} + \lambda_f^+. \label{discrEq:discrAdjEqAdjEndDiffCalcSysEq002}  
\end{eqnarray}

From (\ref{discrEq:discrAdjEqAdjIntDiffCalcSysEq002})-(\ref{discrEq:discrAdjEqAdjIntAlgCalcSysEq002}) we have that  $ l_{f1}(k) = 0, l_{g1}(k) = 0, \ldots,l_{fs}(k) = 0,l_{gs}(k) = 0$ for steps $ k = 0,...,K-1 $. During the analysis we also assume that $ l_{f1}(K) = 0, l_{g1}(K) = 0, \ldots,l_{fs}(K) = 0,l_{gs}(K) = 0$. 
Under that assumption (\ref{discrEq:discrAdjEqRIntDiffCalcSysEq002})-(\ref{discrEq:discrAdjEqRIntAlgCalcSysEq002}) is equivalent to
\begin{eqnarray}
r_{fi} &=& \sum_{j=1}^{s} ha_{ji} f_{xi}^T r_{fj} + g_{xi}^T r_{gi}+  hb_if_{xi}^T \lambda_f^+  \label{discrEq:discrAdjEqRIntDiffCalcSysEq003} \\
0 &=& \sum_{j=1}^{s} ha_{ji}  f_{zi}^T r_{fj} + hb_if_{zi}^T \lambda_f^+  \label{discrEq:discrAdjEqRIntAlgCalcSysEq003}
\end{eqnarray}

Let us now introduce the auxiliary variables 
\begin{eqnarray}
\lambda_{fi} &=& \lambda_f^+ + \sum_{j=1}^{s} \frac{a_{ji}}{b_i} r_{fj},  \label{discrEq:discrAdjEqLamIntDiffCalcSysEq001} \\
\lambda_{gi} &=& \frac{r_{gi}}{hb_i} \label{discrEq:discrAdjEqLamIntAlgCalcSysEq001}
\end{eqnarray}
for $ i = 1,\ldots,s $. (\ref{discrEq:discrAdjEqRIntDiffCalcSysEq003})-(\ref{discrEq:discrAdjEqRIntAlgCalcSysEq003}) can then be transformed to
\begin{eqnarray}
r_{fi} &=& hb_i \left( f_{xi}^T \lambda_{fi} + g_{xi}^T \lambda_{gi} \right), \label{discrEq:discrAdjEqRIntDiffCalcSysEq004} \\
0 &=& hb_i f_{zi}^T \lambda_{fi} . \label{discrEq:discrAdjEqRIntAlgCalcSysEq004} 
\end{eqnarray} 
If we put (\ref{discrEq:discrAdjEqRIntDiffCalcSysEq004}) into (\ref{discrEq:discrAdjEqLamIntDiffCalcSysEq001}) we obtain
\begin{eqnarray}
\lambda_{fi} &=& \lambda_f^+ + h\sum_{j=1}^{s} \frac{a_{ji}b_j}{b_i}  \left( f_{xj}^T \lambda_{fj} + g_{xj}^T \lambda_{gj} \right). \label{discrEq:discrAdjEqLamIntDiffCalcSysEq002} 
\end{eqnarray}
Dividing (\ref{discrEq:discrAdjEqRIntAlgCalcSysEq004}) by $ hb_i $ results in
\begin{eqnarray}
0 &=& f_{zi}^T \lambda_{fi}. \label{discrEq:discrAdjEqLamIntAlgCalcSysEq002}
\end{eqnarray}
If we put (\ref{discrEq:discrAdjEqRIntDiffCalcSysEq004}) into (\ref{discrEq:discrAdjEqAdjEndDiffCalcSysEq002}) we obtain
\begin{eqnarray}
\lambda_f &=& \lambda_f^+ + h\sum_{i=1}^{s} b_i \left( f_{xi}^T \lambda_{fi} + g_{xi}^T \lambda_{gi} \right). \label{discrEq:discrAdjEqAdjEndDiffCalcSysEq003}  
\end{eqnarray}

Let us rewrite equations (\ref{discrEq:discrAdjEqLamIntDiffCalcSysEq002})-(\ref{discrEq:discrAdjEqAdjEndDiffCalcSysEq003}) by introducing the discrete step argument
\begin{eqnarray}
\lambda_{fi}(k) 
&=& \lambda_f(k+1) + h(k) \sum_{j=1}^{s} \frac{a_{ji}b_j}{b_i} \left[  f_x^T 
\left( 
x_j\left( k+1 \right) , z_j\left( k+1 \right), u(k) 
\right) \lambda_{fj}(k) \right.
\nonumber\\
&&\qquad\qquad\qquad\qquad\qquad\qquad\ \
\left.  +g_x^T 
\left( 
x_i\left( k+1 \right) \right) \lambda_{gj}(k) \right], \label{discrEq:discrAdjEqLamIntDiffCalcSysEq004}\\
0 &=& f^T_z\left( x_i\left( k+1 \right) \right) \lambda_{fi}(k)\label{discrEq:discrAdjEqLamIntAlgCalcSysEq004}
\end{eqnarray}
for $ i = 1,\dots,s $ and
\begin{eqnarray}
\lambda_{f}(k)  &=& \lambda_f(k+1) + h(k)\sum_{i=1}^{s} b_i \left[ f_x^T 
\left( 
x_j\left( k+1 \right) , z_i\left( k+1 \right), u(k) 
\right) \lambda_{fi}(k) \right.
\nonumber\\
&&\qquad\qquad\qquad\qquad\qquad\qquad\ \
\left.  +g_x^T 
\left( 
x_i\left( k+1 \right) \right) \lambda_{gi}(k) \right] \label{discrEq:discrAdjEqAdjEndDiffCalcSysEq004} 
\end{eqnarray}
and augment them 
by the equation for $ \lambda_{g}(k) $
\begin{eqnarray}
\lambda_{g}(k) &=& \lambda_{g}(k+1) +\sum_{i=1}^{s} b_i^- (\lambda_{gi}(k) - \lambda_{g}(k+1) ). \label{discrEq:discrAdjEqAdjEndAlgCalcSysEq004} 
\end{eqnarray}

We now consider the Runge-Kutta scheme (\ref{discrEq:rkSchemeIntAdjDiff001})-(\ref{discrEq:rkSchemeStepAdjAlg001}) under the assumption that 
the discrete steps are the same as in the forward scheme and the Runge-Kutta scheme coefficients satisfy $ \bar{a}_{ij} = \frac{a_{ji}b_j}{b_i}, \bar{b}_i = b_i, \bar{c}_i = 1-c_i $  
\begin{eqnarray}
\lambda_{fi}(k) 
&=& \lambda_f(k+1) + h(k) \sum_{j=1}^{s} \frac{a_{ji}b_j}{b_i} \left[  f_x^T 
\left( 
x\left( t_j(k) \right) , z\left( t_j(k) \right), u(k) 
\right) \lambda_{fj}(k) \right.
\nonumber\\
&&\qquad\qquad\qquad\qquad\qquad\qquad\ \
\left.  +g_x^T 
\left( 
x\left( t_j(k) \right) \right) \lambda_{gj}(k) \right], \label{discrEq:rkSchemeIntAdjDiff002}\\
0 &=& -f^T_z\left( x\left( t_i(k) \right) \right) \lambda_{fi}(k) \label{discrEq:rkSchemeIntAdjAlg002}
\end{eqnarray}
for $ i = 1,\dots,s $ and
\begin{eqnarray}
\lambda_{f}(k)  &=& \lambda_f(k+1) + h(k) \sum_{i=1}^{s} b_{i} \left[  f_x^T 
\left( 
x\left( t_i(k) \right) , z\left( t_i(k) \right), u(k)
\right) \lambda_{fi}(k) \right.
\nonumber\\
&&\qquad\qquad\qquad\qquad\qquad\qquad\ \
\left.  +g_x^T 
\left( 
x\left( t_i(k) \right) \right) \lambda_{gi}(k) \right], \label{discrEq:rkSchemeStepAdjDiff002} \\
\lambda_{g}(k) &=& \lambda_{g}(k+1) +\sum_{i=1}^{s} b_i^- (\lambda_{gi}(k) - \lambda_{g}(k+1) ), \label{discrEq:rkSchemeStepAdjAlg002}  
\end{eqnarray}
where
\begin{equation}
t_i(k) = t(k) + c_ih(k).
\end{equation}
Now the Runge-Kutta scheme (\ref{discrEq:rkSchemeIntAdjDiff002})-(\ref{discrEq:rkSchemeStepAdjAlg002}) for adjoint equations is almost identical to discrete adjoint equations (\ref{discrEq:discrAdjEqLamIntDiffCalcSysEq004})-(\ref{discrEq:discrAdjEqAdjEndAlgCalcSysEq004}). The only difference is that in discrete equations $ x_i(k+1) $ and $ z_i(k+1) $ are used instead of $ x(t(k)+c_ih(k)) $ and $ z(t(k)+c_ih(k)) $. The convergence of (\ref{discrEq:rkSchemeIntAdjDiff002})-(\ref{discrEq:rkSchemeStepAdjAlg002}) is guaranteed if only the Runge-Kutta scheme $ \bar{a}_{ij} ,\ \bar{b}_i,\ \bar{c}_i,\ i,j = 1,\dots,s $ satisfies appropriate conditions. We will also justify that the usage of $ x_i(k+1) $ and $ z_i(k+1) $ instead of $ x(t(k)+c_ih(k)) $ and $ z(t(k)+c_ih(k)) $ does not destroy the convergence of (\ref{discrEq:discrAdjEqLamIntDiffCalcSysEq004})-(\ref{discrEq:discrAdjEqAdjEndAlgCalcSysEq004}) solutions to continuous adjoint trajectory. 

The partial derivative $ F_{u}(k)  $ is
\begin{equation}
F_{u}(k) = \left(
\begin{array}{c}
- h\sum_{j=1}^{s}a_{1j} f_{uj}  \\
0 \\
\vdots  \\
- h\sum_{j=1}^{s}a_{sj} f_{uj} \\
0 \\
- h\sum_{i=1}^{s}b_{i} f_{ui}
\end{array}
\right ),
\end{equation}
where for the sake of the shorter notation we omitted the discrete step argument and introduced
\begin{equation}
f_{ui} = f_u(x_i,u).
\end{equation}
Let us now inspect $ -F_u^T(k)R(k) $:  
\begin{eqnarray}
-F_u^T(k)R(k) &=&  
- \left(
\begin{array}{c}
- h\sum_{j=1}^{s}a_{1j} f_{uj}  \\
0 \\
\vdots  \\
- h\sum_{j=1}^{s}a_{sj} f_{uj} \\
0 \\
- h\sum_{i=1}^{s}b_{i} f_{ui}
\end{array}
\right )^T
\left( 
\begin{array}{c}
r_{f1} \\ r_{g1} \\ \vdots \\ 
r_{fs} \\ r_{gs} \\ r_f
\end{array}
\right) \\
&=&  \sum_{i = 1}^{s}  h \left(\sum_{j=1}^{s}a_{ij}  f^T_{uj} \right) r_{fi} + h\sum_{i=1}^{s}b_{i} f^T_{ui} r_f. \nonumber
\end{eqnarray}
From (\ref{discrEq:discrAdjEqREndDiffCalcSysEq001}) and (\ref{discrEq:discrAdjEqRIntDiffCalcSysEq004}) we get
\begin{eqnarray}
&&-F_u^T(k)R(k) =  \\
&&\sum_{i = 1}^{s}  h \left(\sum_{j=1}^{s}a_{ij}  f^T_{uj} \right) hb_i \left( f_{xi}^T \lambda_{fi} + g_{xi}^T \lambda_{gi} \right) + h\sum_{i=1}^{s}b_{i} f^T_{ui} \lambda_f^+ = \nonumber \\
&& h \sum_{j = 1}^{s} f^T_{uj} \left( h \sum_{i=1}^{s}a_{ij} b_i \left( f_{xi}^T \lambda_{fi} + g_{xi}^T \lambda_{gi} \right) \right)  + h\sum_{i=1}^{s}b_{i} f^T_{ui} \lambda_f^+ = \nonumber \\
&& h \sum_{j = 1}^{s} f^T_{uj}b_j \left( h \sum_{i=1}^{s} \frac{a_{ij} b_i}{b_j} \left( f_{xi}^T \lambda_{fi} + g_{xi}^T \lambda_{gi} \right)\right)  + h\sum_{i=1}^{s}b_{i} f^T_{ui} \lambda_f^+ . \nonumber
\end{eqnarray}
From (\ref{discrEq:discrAdjEqLamIntDiffCalcSysEq002}) we get that
\begin{eqnarray}
-F_u^T(k)R(k) &=&  h \sum_{j = 1}^{s} f^T_{uj}b_j \left( \lambda_{fj} - \lambda^+_f \right)  + h\sum_{i=1}^{s}b_{i} f^T_{ui} \lambda_f^+ \\
&=&  h \sum_{i = 1}^{s} f^T_{ui}b_i \left(  \lambda_{fi} - \lambda^+_f \right)  + h\sum_{i=1}^{s}b_{i} f^T_{ui} \lambda_f^+ \nonumber\\
&=&  h \sum_{i = 1}^{s} b_i f^T_{ui}  \lambda_{fi} . \nonumber
\end{eqnarray} 
The discrete reduced gradient now takes the form
\begin{equation}
\frac{\phi(X(K))}{du_n} =  \sum_{k=k_{n-1}}^{k_n-1} h(k) \sum_{i = 1}^{s} b_i f_u^T(x_i(k+1),u_n)  \lambda_{fi}(k).
\label{discrEq:discrRedGradDAE001}
\end{equation}

Under the assumption that the discrete steps are the same as in the forward scheme and the quadrature coefficients satisfy $  \tilde{b}_i = b_i, \tilde{c}_i = c_i$, the reduced gradients formula (\ref{discrEq:redGradQuadContSlid}) is
\begin{equation}
\frac{d \phi(x(t_f))}{ du_n } \simeq \sum_{k=k_{n-1}}^{k_n-1} h(k) \sum_{i=1}^{s} b_i f_u^T \left(x \left(t(k)+c_i h(k)\right),u_n \right) \lambda_f \left( t(k)+c_i h(k) \right).
\label{discrEq:redGradQuadContSlid001}
\end{equation}
Now the quadrature scheme (\ref{discrEq:redGradQuadContSlid001}) for continuous reduced gradients is almost identical to discrete reduced gradients  (\ref{discrEq:discrRedGradDAE001}). The only difference is that in (\ref{discrEq:discrRedGradDAE001}) $ x_i(k+1) $ and $ \lambda_{fi}(k) $ are used instead of $ x(t(k)+c_ih(k)) $ and $ \lambda_f(t(k)+c_ih(k)) $. The convergence of (\ref{discrEq:redGradQuadContSlid001}) is guaranteed if only the quadrature with coefficients $ b_i,\ c_i,\ i = 1,\dots,s $ satisfies appropriate conditions. We will also justify that the usage of  $ x_i(k+1) $ and $ \lambda_{fi}(k) $ instead of $ x(t(k)+c_ih(k)) $ and $ \lambda_f(t(k)+c_ih(k)) $ does not destroy the convergence of discrete reduced gradients (\ref{discrEq:discrRedGradDAE001}) to continuous reduced gradients. 

Now we want to justify the convergence of discrete adjoint trajectories and discrete reduced gradients to their continuous counterparts. To obtain the right error order estimated we always assume that the system functions $ f(x,z,u) $ and $ g(x) $ and their partial derivatives are Lipschitz continuous functions. To achieve that goal we use the theorems presented in \cite{hlr89}. That theorems are formulated under the constant step size assumption,  
\begin{equation}
h(k) = h,\ k=1,\ldots,K ,
\end{equation}
and we also carry out our analysis under that assumption. The variable step size case can be cumbersome to analyze, and in this section we will make a short note on that problem. 

We assume that the system equations are integrated using RADAU IIA scheme for $ s=3 $. The coefficients that appear in discrete adjoint equations $ \bar{a}_{ij} = \frac{a_{ji}b_j}{b_i}, \bar{b}_i = b_i, \bar{c}_i = 1-c_i $ define  the RADAU IA scheme (see \cite{ps19a}). The global error of the differential state $ x(t) $ is (\cite{hlr89}, Theorem 4.4 p. 36 and Theorem 5.9 p. 67)
\begin{equation}
x(k) - x(t(k)) = O(h^p). \label{discrEq:errEstStateDiffStep} 
\end{equation}
and the global error of the algebraic state $ z(t) $ is (\cite{hlr89}, Theorem 4.6 p. 40)
\begin{equation}
z(k) - z(t(k)) = O(h^q). \label{discrEq:errEstStateAlgStep} 
\end{equation}
where for the RADAU IIA scheme and $ s=3 $ we have \cite{hlr89}
\begin{equation}
p=5,\  q=3.
\end{equation}
We emphasize that in the DAEs case the invertibility of the matrix of coefficients $ A = (a_{ij}) $ is essential for convergence of the scheme. Also the value of so called radius of stability 
\begin{equation}
R(\infty) = 1-b^TA^{-1}\mathbf{1}.
\end{equation}
plays an important role. For RADAU IA and RADAU IIA schemes $ A $ matrices are invertible and $ R(\infty)=0 $, so the assumptions of the appropriate theorems from \cite{hlr89} are satisfied. 

To derive the subsequent results we consider the unperturbed Runge -- Kutta scheme
\begin{eqnarray}
x_i^n(k+1) &=& x^n(k)+h\sum_{j=1}^{s}a_{ij}f \left( x_j^n(k+1),\ z_j^n(k+1), u(k) \right), \label{discrEq:rkSchemeIntDiffUnpert001} \\
0 &=& g \left( x_i^n(k+1)  \right), \label{discrEq:rkSchemeIntAlgUnpert001} 
\end{eqnarray}
for $ i = 1,...,s $ and 
\begin{eqnarray}
x^n(k+1) &=& x^n(k) +h\sum_{i=1}^{s}b_i f \left( x_i^n(k+1),\ z_i^n(k+1), u(k) \right) \label{discrEq:rkSchemeStepDiffUnpert001} 
\\ z^n(k+1) &=& z^n(k) +\sum_{i=1}^{s}  b_i^-  (z_i^n(k+1) - z^n(k)). \label{discrEq:rkSchemeStepAlgUnpert001} 
\end{eqnarray} 
and the perturbed Runge-Kutta scheme
\begin{eqnarray}
x_i^p(k+1) &=& x^n(k)+h\sum_{j=1}^{s}a_{ij}f \left( x_j^p(k+1),\ z_j^p(k+1), u(k) \right) +\delta_i, \label{discrEq:rkSchemeIntDiffPert001} \\
0 &=& g \left( x_i^p(k+1)  \right), \label{discrEq:rkSchemeIntAlgPert001} 
\end{eqnarray}
for $ i = 1,...,s $ and 
\begin{eqnarray}
x^p(k+1) &=& x^n(k) +h\sum_{i=1}^{s}b_i f \left( x_i^p(k+1),\ z_i^p(k+1), u(k) \right)+\delta_{s+1}. \label{discrEq:rkSchemeStepDiffPert001} 
\\ z^p(k+1) &=& z^n(k) +\sum_{i=1}^{s}  b_i^-  (z_i^p(k+1) - z^n(k))+\delta_{s+2}. \label{discrEq:rkSchemeStepAlgPert001} 
\end{eqnarray} 
From (\cite{hlr89}, Theorem 4.2 p. 33) we have
\begin{eqnarray}
x_i^p(k+1) - x_i^n(k+1) &=& O(\delta) \label{discrEq:errEstStateDiffIntPert}\\
z_i^p(k+1) - z_i^n(k+1) &=& \frac{1}{h}O(\delta) \label{discrEq:errEstStateAlgIntPert}
\end{eqnarray}
where $ \delta = max\{\delta_1,\dots,\delta_s \} $. Contrary to the ODEs case, using the variable step size for  integration of DAEs may be cumbersome, because of the factor $ \frac{1}{h} $ that occurs in (\ref{discrEq:errEstStateAlgIntPert}). The estimate (\ref{discrEq:errEstStateAlgIntPert}) indicates that using very small step sizes may lead to big errors. On the other hand, the results for the constant step size can be extended to the variable step size case if we assume that the following condition holds
\begin{equation}
h_m \geq h(k) \geq \mu h_m,\ k=1,\ldots,K
\end{equation}   
where $ h_m $ is the maximum step size and $ \mu\in(0,1) $ is a constant independent from $ h_m $. 

The main step errors can be estimated as
\begin{eqnarray}
x^p(k+1) - x^n(k+1) &=& O(\delta)+O(\delta_{s+1}) \label{discrEq:errEstStateDiffStepPert}\\
z^p(k+1) - z^n(k+1) &=& \frac{1}{h}O(\delta)+O(\delta_{s+2}) \label{discrEq:errEstStateAlgStepPert}
\end{eqnarray}

Let us now derive the estimates of the global errors $  x_i(k+1) - x(t(k) + c_ih(k))  $ and $  z_i(k+1) - z(t(k) + c_ih(k))  $. The nominal Runge-Kutta scheme is formulated assuming that the exact state value is known at $ t(k) $ 
\begin{eqnarray}
x_i^n(k+1) &=& x(t(k))+h\sum_{j=1}^{s}a_{ij}f \left( x_j^n(k+1),\ z_j^n(k+1), u(k) \right), \label{discrEq:rkSchemeIntDiffUnpert002} \\
0 &=& g \left( x_i^n(k+1)  \right), \label{discrEq:rkSchemeIntAlgUnpert002} 
\end{eqnarray}
for $ i = 1,...,s $. The following local error estimates are valid for RADAU IIA scheme (\cite{hlr89}, Lemma 4.3 p. 34) 
\begin{eqnarray}
x^n_i(k+1) - x(t(k) + c_ih)  &=& O(h^{q+1}) \label{discrEq:errEstStateDiffIntLoc} \\
z^n_i(k+1) - z(t(k) + c_ih)  &=& O(h^{q}) \label{discrEq:errEstStateAlgIntLoc}
\end{eqnarray}
In this case the perturbed Runge-Kutta scheme is actually the regular Runge-Kutta scheme, with the approximation of the state $ x(k) $ at time $ t(k) $ used
\begin{eqnarray}
&& x_i(k+1) = \nonumber\\
&& x(k)+h\sum_{j=1}^{s}a_{ij}f \left( x_j(k+1),\ z_j(k+1), u(k) \right) = \label{discrEq:rkSchemeIntDiffPert002} \\
&& x(t(k))+h\sum_{j=1}^{s}a_{ij}f \left( x_j(k+1),\ z_j(k+1), u(k) \right) + (x(k)-x(t(k)), \nonumber \\
&& 0 = g \left( x_i(k+1)  \right), \label{discrEq:rkSchemeIntAlgPert002} 
\end{eqnarray}
for $ i = 1,...,s $, so the perturbations are
\begin{equation}
\delta_i = \delta = x(k)- x(t(k)),\ i = 1,\dots,s.
\end{equation}
From (\ref{discrEq:errEstStateDiffStep}) we have $  \delta  =  x(k)- x(t(k))  =  O(h^p) $. From (\ref{discrEq:errEstStateDiffIntPert}) and (\ref{discrEq:errEstStateAlgIntPert}) we obtain
\begin{eqnarray}
x_i(k+1) - x^n_i(k+1) &=& O(h^p). \label{discrEq:errEstStateDiffIntPert001} \\
z_i(k+1) - z^n_i(k+1) &=& O(h^{p-1}). \label{discrEq:errEstStateAlgIntPert001}
\end{eqnarray} 
By combining (\ref{discrEq:errEstStateDiffIntLoc})-(\ref{discrEq:errEstStateAlgIntLoc}) and (\ref{discrEq:errEstStateDiffIntPert001})-(\ref{discrEq:errEstStateAlgIntPert001})  we obtain the required global error estimates
\begin{eqnarray}
x_i(k+1) - x(t(k) + c_ih) &=& O(h^{q+1})+ O(h^p) = O(h^{min\{p,q+1\}})\nonumber \\
& = & O(h^{q_x^g}) \label{discrEq:errEstDiffIntGlob} \\
z_i(k+1) - z(t(k) + c_ih) &=& O(h^{q})+ O(h^{p-1})\nonumber \\
& = & O(h^{q_z^g}) \label{discrEq:errEstAlgIntGlob} 
\end{eqnarray}
For the RADAU IIA scheme with $ s=3 $ we get
\begin{equation}
q_x^g = min\{5,3+1\} = 4,\ q_z^g = min\{5-1,3\} = 3.\label{discrEq:errEstDiffIntGlobR2A} \\
\end{equation}

Let us now define the unperturbed Runge--Kutta scheme as the Runge-Kutta scheme for continuous adjoint equations with the exact adjoint state known at time $ t(k+1) $
\begin{eqnarray}
\lambda_{fi}^n(k) 
&=& \lambda_f(t(k+1)) + h \sum_{j=1}^{s} \bar{a}_{ij} \left[  f_x^T 
\left( 
x\left( t_j(k) \right) , z\left( t_j(k) \right), u(k) 
\right) \lambda_{fj}^n(k) \right.
\nonumber\\
&&\qquad\qquad\qquad\qquad\qquad\qquad\ \
\left.  +g_x^T 
\left( 
x\left( t_j(k) \right) \right) \lambda_{gj}^n(k) \right], \label{discrEq:rkSchemeIntAdjDiff003}\\
0 &=& -f^T_z\left( x\left( t_i(k) \right) \right) \lambda_{fi}^n(k) \label{discrEq:rkSchemeIntAdjAlg003}
\end{eqnarray}
for $ i = 1,\dots,s $ and
\begin{eqnarray}
\lambda_{f}^n(k)  &=& \lambda_f(t(k+1)) + h \sum_{i=1}^{s} \bar{b}_{i} \left[  f_x^T 
\left( 
x\left( t_i(k) \right) , z\left( t_i(k) \right), u(k)
\right) \lambda_{fi}^n(k) \right.
\nonumber\\
&&\qquad\qquad\qquad\qquad\qquad\qquad\ \
\left.  +g_x^T 
\left( 
x\left( t_i(k) \right) \right) \lambda_{gi}^n(k) \right], \label{discrEq:rkSchemeStepAdjDiff003} \\
\lambda_{g}^n(k) &=& \lambda_{g}(t(k+1)) +\sum_{i=1}^{s} \bar{b}_i^- (\lambda_{gi}^n(k) - \lambda_{g}(t(k+1)) ), \label{discrEq:rkSchemeStepAdjAlg003}  
\end{eqnarray}
The coefficients $ \bar{a}_{ij}, \bar{b}_{i}, \bar{c}_{i}$ constitutes RADAU IA scheme (\cite{ps19a}), so the following local error estimates are valid (\cite{hlr89}, Theorem 4.3 p. 43)
\begin{eqnarray}
\lambda^n_{fi}(k) - \lambda_f(t(k)+c_ih) &=& O(h^{\bar{q}+1}),\ i=1,\dots,s\ , \label{discrEq:errEstDiffAdjIntLoc} \\
\lambda^n_{gi}(k) - \lambda_g(t(k)+c_ih) &=& O(h^{\bar{q}}),\ i=1,\dots,s\ , \label{discrEq:errEstAlgAdjIntLoc} \\
\lambda^n_{f}(k) - \lambda_f(t(k)) &=& O(h^{\bar{q}+1}), \label{discrEq:errEstDiffAdjEndLoc} \\
\lambda^n_{g}(k) - \lambda_g(t(k)) &=& O(h^{\bar{q}}). \label{discrEq:errEstAlgAdjEndLoc} 
\end{eqnarray}
where for the RADAU IA scheme with $ s=3 $ we have (\cite{hlr89})
\begin{equation}
\bar{p}=5,\ \bar{q}=2.   
\end{equation}

The perturbed Runge--Kutta scheme is defined as the discrete adjoint equations with the exact adjoint state known at time $ t(k+1) $ 
\begin{eqnarray}
\lambda_{fi}^p(k) 
&=& \lambda_f(t(k+1)) + h \sum_{j=1}^{s} \bar{a}_{ij} \left[  f_x^T 
\left( 
x_j(k+1) , z_j(k+1), u(k) 
\right) \lambda_{fj}^p(k) \right.
\nonumber\\
&&\qquad\qquad\qquad\qquad\qquad\qquad\ \
\left.  +g_x^T 
\left( 
x_j(k+1) \right) \lambda_{gj}^p(k) \right], \label{discrEq:discrAdjEqIntDiffCalcSysEq004}\\
0 &=& -f^T_z\left( x_i(k+1) \right) \lambda_{fi}^p(k) \label{discrEq:discrAdjEqIntAlgCalcSysEq004}
\end{eqnarray}
for $ i = 1,\dots,s $ and
\begin{eqnarray}
\lambda_{f}^p(k)  &=& \lambda_f(t(k+1)) + h \sum_{i=1}^{s} \bar{b}_{i} \left[  f_x^T 
\left( 
x_i(k+1) , z_i(k+1), u(k)
\right) \lambda_{fi}^p(k) \right.
\nonumber\\
&&\qquad\qquad\qquad\qquad\qquad\qquad\ \
\left.  +g_x^T 
\left( 
x_i(k+1) \right) \lambda_{gi}^p(k) \right], \label{discrEq:discrAdjEqAdjEndDiffCalcSysEq005} \\
\lambda_{g}^p(k) &=& \lambda_{g}(t(k+1)) +\sum_{i=1}^{s} \bar{b}_i^- (\lambda_{gi}^p(k) - \lambda_{g}(t(k+1)) ), \label{discrEq:discrAdjEqAdjEndAlgCalcSysEq005}  
\end{eqnarray}
The perturbed system can be rewritten as
\begin{eqnarray}
\lambda_{fi}^p(k) 
&=& \lambda_f(t(k+1)) + h \sum_{j=1}^{s} \bar{a}_{ij} \left[  f_x^T 
\left( 
x(t_j(k)) , z(t_j(k)), u(k) 
\right) \lambda_{fj}^p(k) \right.
\nonumber\\
&&\qquad\qquad\qquad\qquad\qquad\qquad\ \
\left.  +g_x^T 
\left( 
x(t_j(k)) \right) \lambda_{gj}^p(k) \right] + \delta_i, \label{discrEq:discrAdjEqIntDiffCalcSysEq005}\\
0 &=& -f^T_z\left( x(t_i(k)) \right) \lambda_{fi}^p(k) \label{discrEq:discrAdjEqIntAlgCalcSysEq005}
\end{eqnarray}
for $ i = 1,\dots,s $ and
\begin{eqnarray}
\lambda_{f}^p(k)  &=& \lambda_f(t(k+1)) + h \sum_{i=1}^{s} \bar{b}_{i} \left[  f_x^T 
\left( 
x(t_i(k)) , z(t_i(k)), u(k)
\right) \lambda_{fi}^p(k) \right.
\nonumber\\
&&\qquad\qquad\qquad\qquad\qquad\qquad\ \
\left.  +g_x^T 
\left( 
x(t_i(k)) \right) \lambda_{gi}^p(k) \right] +\delta_{s+1}, \label{discrEq:discrAdjEqAdjEndDiffCalcSysEq006} \\
\lambda_{g}^p(k) &=& \lambda_{g}(t(k+1)) +\sum_{i=1}^{s} \bar{b}_i^- (\lambda_{gi}^p(k) - \lambda_{g}(t(k+1)) ), \label{discrEq:discrAdjEqAdjEndAlgCalcSysEq006}  
\end{eqnarray}
where the perturbations result from differences $ x_i(k+1) - x\left( t(k)+c_ih(k) \right) = O(h^{q_x^g}) = O(h^4)$ and $ z_i(k+1) - z\left( t(k)+c_ih(k) \right) = O(h^{q_z^g}) = O(h^3)$. The perturbations satisfy 
\begin{equation}
\delta_i = O(h^{q_z^g}),\ i=1,\dots,s+1.
\end{equation}
From (\ref{discrEq:errEstStateDiffIntPert})-(\ref{discrEq:errEstStateAlgStepPert}) we obtain the error estimates
\begin{eqnarray}
\lambda^p_{fi}(k) - \lambda^n_{fi}(k) &=& O(h^{q_z^g}), \label{discrEq:errEstDiscrAdjDiffIntPert} \\
\lambda^p_{gi}(k) - \lambda^n_{gi}(k) &=& O(h^{q_z^g-1}), \label{discrEq:errEstDiscrAdjAlgIntPert} \\
\lambda^p_f(k) - \lambda^n_f(k) &=& O(h^{q_z^g}), \label{discrEq:errEstDiscrAdjDiffEndPert} \\
\lambda^p_g(k) - \lambda^n_g(k) &=& O(h^{q_z^g-1}). \label{discrEq:errEstDiscrAdjAlgEndPert}
\end{eqnarray}
By combining (\ref{discrEq:errEstDiffAdjIntLoc})-(\ref{discrEq:errEstAlgAdjEndLoc}) with (\ref{discrEq:errEstDiscrAdjDiffIntPert})-(\ref{discrEq:errEstDiscrAdjAlgEndPert}) we obtain the following local error estimates
\begin{eqnarray}
\lambda^p_{fi}(k) - \lambda_f(t(k)+c_ih) &=& O(h^{min\{q_z^g,\bar{q}+1\}})=O(h^{\bar{q}_f^l}),\ i=1,\dots,s\ , \label{discrEq:errEstDiscrAdjDiffIntLoc} \\
\lambda^p_{gi}(k) - \lambda_g(t(k)+c_ih) &=& O(h^{min\{q_z^g-1,\bar{q}\}})=O(h^{\bar{q}_g^l}),\ i=1,\dots,s\ , \label{discrEq:errEstDiscrAdjAlgIntLoc} \\
\lambda^p_{f}(k) - \lambda_f(t(k)) &=& O(h^{min\{q_z^g,\bar{q}+1\}})=O(h^{\bar{p}_f^l}), \label{discrEq:errEstDiscrAdjDiffEndLoc} \\
\lambda^p_{g}(k) - \lambda_g(t(k)) &=& O(h^{min\{q_z^g-1,\bar{q}\}})=O(h^{\bar{p}_g^l}). \label{discrEq:errEstDiscrAdjAlgEndLoc} 
\end{eqnarray}
For the RADAU IA scheme with $ s=3 $ we have 
\begin{eqnarray}
\bar{q}_f^l = \bar{p}_f^l &=& min\{q_z^g,\bar{q}+1\} = min\{3,2+1\} =3 \\
\bar{q}_g^l = \bar{p}_g^l &=& min\{q_z^g-1,\bar{q}\} = min\{3-1,2\} =2 
\end{eqnarray}
The global error of $ \lambda^p_{f}(k) $ can be obtained from (\cite{hlr89}, Theorem 4.4 p. 36)
\begin{equation}
\lambda_{f}(k) - \lambda_f(t(k)) = O(h^{\bar{p}_f^l-1}) = O(h^{\bar{p}_f^g}). \label{discrEq:errEstDiscrAdjDiffStepGlob}
\end{equation} 
For RADAU IIA scheme with $ s=3 $ we obtain $ \bar{p}_f^g = min\{q_z^g,\bar{q}+1\}-1 = min\{3,2+1\}-1 = 2$.
Let us now derive the global error estimate of $ \lambda_{fi}(k) - \lambda_f(t(k)+c_ih(k)) $. Using the global error estimate (\ref{discrEq:errEstDiscrAdjDiffStepGlob}) and the local error estimate (\ref{discrEq:errEstDiscrAdjDiffIntPert}) we can repeat the reasoning for discrete system equations and derive the global error estimate
\begin{equation}
\lambda_{fi}(k) - \lambda_f(t(k)+c_ih(k)) = O(h^{min\{\bar{p}_f^g, \bar{q}_f^l\} }) = O(h^{ \bar{q}_f^g}). \label{discrEq:errEstDiscrAdjDiffIntGlob}
\end{equation}
For RADAU IIA scheme with $ s=3 $ we obtain $ \bar{q}_f^g = $ $min\{\bar{p}_f^g, \bar{q}_f^l\} $ $= min\{2,3\} = 2$.

Having the global error estimates (\ref{discrEq:errEstDiffIntGlob}) for $ x_i(k+1) - x(t_i(k))$ and (\ref{discrEq:errEstDiscrAdjDiffIntGlob}) for $ \lambda_{fi}(k) - \lambda_{f}(t_i(k)) $, we can repeat the derivation of the reduced gradients calculation order presented in the first part of the paper (for ODEs case) to obtain
\begin{equation}
\frac{d \phi(X(K))}{du_n} - \frac{d \phi(x(t_f))}{ du_n } =  O(h^{\tilde{p}_d}) 
\end{equation}  
where $ \tilde{p}_d = min\{\tilde{p}, q_x^g, \bar{q}_f^g \} $. For RADAU IIA scheme we obtain
\begin{equation}
\tilde{p}_d = min\{5, 4,2 \} =  2.
\end{equation} 
This result confirms that the discrete reduced gradients converge to the continuous reduced gradients with an order at least $ \tilde{p}_d =2 $ (for RADAU IIA scheme with $ s=3 $). The discrete reduced gradients provide therefore an efficient and reliable method for approximation of continuous reduced gradients. 

In Table \ref{Orders} we summarize the proven minimal integration orders for system and adjoint equations derived in this paper and in \cite{ps19a}. It follows that if a state trajectory includes sections with sliding modes then the reduced gradient can be approximated with the accuracy at least $O(h^2)$, otherwise the accuracy is at least $O(h^3)$. This accuracy estimate takes into account the fact that if a discrete state changes then we have to start the integration of system and adjoint equations with the perturbed initial values of system and adjoint variables respectively---the influence of that perturbation on the accuracy of system and adjoint variables determination follows from Theorem 4.3 in \cite{hw96} (cf. (120)--(121) and (139) in \cite{ps19a}), or from Theorem 4.2 in \cite{hlr89} (cf. (\ref{discrEq:errEstStateDiffStepPert})--(\ref{discrEq:errEstStateAlgStepPert}) and  (\ref{discrEq:errEstDiscrAdjDiffIntPert})--(\ref{discrEq:errEstDiscrAdjAlgEndPert})).

\begin{table}
	\label{Orders}
	\caption{Integration orders of convergence if RADAU IIA is applied to systems equations}
	\begin{center}
		\begin{tabular}{llll}
			\hline{\rule{0pt}{12pt}}
			Equations & Orders $(x,z)$ & Orders $(\lambda_f,\lambda_g)$ & gradient orders \\[2pt]
			\hline{\rule{0pt}{12pt}}
			$x' = f(x,u)$  &     $p=5$ &  $\bar{q}_f = 4$ & $p_d = 3$   \\[2pt] \hline{\rule{0pt}{12pt}} 
			$x' = f(x,u) + g_x(x)z$ & & & \\
			$0 = g(x)$  &   $p = 5$, $q = 3$ & $\bar{q}_f = 2$, $\bar{q}_g = 2$   &  $p_d=2$ \\[2pt]
			\hline
		\end{tabular}
	\end{center}
\end{table}

\section{Calculating adjoint variables jumps} 
\label{secJump}

In this section we want to discuss the correspondence between jump conditions for discrete and continuous adjoint equations at transition times. Let us consider the hybrid system trajectory described by (\ref{eqAdjEqsFirstCaseOde})-(\ref{eqAdjEqsFirstCaseDaeAlg}). As a result of the numerical integration we obtain the discrete state equations (\ref{discrEq:discrStateEq}). 

In Section \ref{secIntro} we have stated that at a transition time $t_t$ adjoint variables undergo jumps, the extent of the jump depends on the sequence of discrete variables before and after the jump. The analysis which follows concerns the case of the equations (\ref{eqAdjEqsFirstCaseOde})-(\ref{eqAdjEqsFirstCaseDaeAlg}), in that case the extent of the jump can be determined according to the equations (\ref{odeAdjTransHamilConta})--(\ref{odeAdjTransHamilCont})

These equations can be solved with respect to $\pi_t$ and $\lambda_f(t_t^-)$ giving:
\begin{eqnarray}
&{\displaystyle \pi_t  =  -\frac{ \lambda(t_t^+)^T\left (f^2(x(t_t^+),z(t_t^+),u(t_t^+))-f^1(x(t_t^-),u(t_t^-)) \right )}{g_x(x(t_t)) f^1(x(t_t^-),u(t_t^-))}} \nonumber
\end{eqnarray}
(we have taken into account that $g(x(t_t^+))=0$ ).

Between transitions, a system of differential-algebraic equations is integrated with the help of an appropriate numerical integration scheme. The numerical integration scheme is represented by the equation
\begin{equation}
F(X(k+1),X(k),u(k),h(k))=0
\label{discrTimeEq}
\end{equation}
in which the discrete step $ k $ corresponds to a time instant $ t(k) $.

During the numerical integration of a hybrid trajectory, a possible violation of invariant set conditions has to be monitored. This task is realized by checking the sign changes of 
\begin{equation}
g(x(k))\nonumber
\end{equation} 
in subsequent steps, where $x(k) = x(t(k))$.
When a sign change of $ g(x(k))$ between discrete steps $ k $ and $ k+1 $ is detected, the following problem is solved
\begin{equation}
{\rm find}\ \ t_t\in [t(k),t(k+1)],\ {\rm s.t.}\ \ \hat{g}(t_t) = 0
\label{findTransTime}
\end{equation}
where $ \hat{g}(\cdot) $ is a function, which approximates $g(\cdot))$ on a time interval $ [t(k),t(k+1)] $. When a transition time $ t_t $ is found, the actual iteration of numerical integration is repeated  but with a step-size $ h(k) = t_t-t(k) $ instead of $ h(k) = t(k+1)-t(k) $. The discrete step at which the transition takes place we denote by $k_t$.

We assume that at discrete times $0,\ldots,k_t$ the system evolves according to the equations:
\begin{equation}
F^1(X^1(k+1),X^1(k),u(k),h(k))=0,
\label{discrTimeEq1}
\end{equation}
and at times $k_t,\ldots,K-1$ by the equations
\begin{equation}
F^2(X^2(k+1),X^2(k),u(k),h(k))=0
\label{discrTimeEq2}
\end{equation}

The optimal control problem with that system was investigated in \cite{ps17}. Therein, the adjoint equations for the functional $\phi(x(t_f))$ were established, herein parts of these equations are presented to expose jumps in adjoint variables.  
\begin{subequations}
	\label{adjointEquations}    
	\begin{align}
	{\rm for}\ k = N-1,\ldots,k_{t}+1& \nonumber\\ 
	\Lambda^2(k)&=- F^2_X(k_t)^T\left [F^2_{X^+}(k_t)\right
	]^{-T}\Lambda^2(k+1)
    \label{adjointEquationsStep1}\\
	\Lambda^{2+}(k_t)&=- F^2_X(k_t)^T\left [F^2_{X^+}(k_t)\right
	]^{-T}\Lambda^2(k_t+1)
	\label{adjointEquationsStepPlus}\\
	\tilde{\Lambda}^{1-}(k_t)+\pi(k_t) \left ( g_x(k_t) \right )^T &= 
	\tilde{\Lambda}^{2+}(k_t) \label{adjointEquationsTrans1}\\
    \Lambda^{1-}(k_t)^T \left [ F^1_{X^+}(k_t-1)\right ]^{-1}F^1_h(k_t-1) & = \Lambda^2(k_t+1)^T\left [F^2_{X^+}(k_t)\right ]^{-1}F^2_h(k_t),
    \label{adjointEquationsTrans2}\\
	\Lambda^1(k_t-1)&=- F^1_X(k_t)^T\left [F^1_{X^+}(k_t)\right
	]^{-T}\Lambda^{1-}(k_t)
	\label{adjointEquationsStepMinus}\\
	{\rm for}\ k = k_t-2,\ldots,1 \nonumber \\ 
	\Lambda^1(k)&=- F^1_X(k_t)^T\left [F^1_{X^+}(k_t)\right
	]^{-T}\Lambda^1(k+1).
	\label{adjointEquationsStep2}
	\end{align}
\end{subequations} 
Here, $\tilde{\Lambda}^{1-}(k_t)$ and $\tilde{\Lambda}^{2+}(k_t)$  are parts of vectors $\Lambda^1(k_t)$ and $\Lambda^2(k_t)$ respectively, defined in such a way to be able to extract the essential part of Eq. (32d) in \cite{ps17}. Notice that 
\begin{eqnarray}
\Lambda^2(k) & = & \left( \left (l_{f1}^2(k)\right )^T,\left (l_{g1}^2(k)\right )^T, \ldots, \left (l_{fs}^2(k)\right )^T,\left (l_{gs}^2(k)\right )^T, \left (\lambda_f^2(k)\right )^T \right)^T,\nonumber 
\end{eqnarray}
for $k=k_t+1,\ldots,K$, but 
\begin{eqnarray}
l_{fi}^2(k) &=& 0 \label{discrEq:discrAdjEqAdjIntDiffCalcSysEq001}\\
l_{gi}^2(k) &=& 0 
\label{discrEq:discrAdjEqAdjIntAlgCalcSysEq001}
\end{eqnarray}
for $ i = 1,\ldots,s $ and $k=k_t+1,\ldots,K$ (the justification of that is given in Section \ref{secDiscrete}), and similarly
\begin{eqnarray}
\Lambda^1(k) & = & \left( \left (l_{1}^1(k)\right )^T, \ldots, \left (l_{s}^1(k)\right )^T, \left (\lambda^1(k)\right )^T \right)^T,\nonumber 
\end{eqnarray}
for $k=0,\ldots,k_t-1$ 
and 
\begin{eqnarray}
l_{i}^1(k) & = & 0 \label{discreteAdjSubVector}
\end{eqnarray}
for $ i=1,\ldots,s$ and steps $ k = 0,...,k_t$ (the justification for that is given by Eq. (78) of Section 5 in \cite{ps19a}).

Therefore, we can take 
\begin{eqnarray}
\tilde{\Lambda}^{1-}(k_t) & = & \lambda^{1-}(k_t) \nonumber \\
\tilde{\Lambda}^{2+}(k_t) & = & \lambda_f^{2+}(k_t) \nonumber 
\end{eqnarray}
and then (\ref{adjointEquationsTrans1}) becomes
\begin{eqnarray}
	\lambda^{1-}(k_t)+\pi(k_t) \left ( g_x(k_t) \right )^T &= & \lambda_f^{2+}(k_t) \label{adjointEquationsTrans1a}
\end{eqnarray}

Our aim is to determine $\pi(k_t)$ on the basis of equations 
(\ref{adjointEquationsTrans1a}) and (\ref{adjointEquationsTrans2}). To this end we need analytical formula for $F^1_X(k_t)$, $F^1_{X^+}(k_t)$, $F^2_X(k_t)$, $F^2_{X^+}(k_t)$, $F^1_h(k_t)$, $F^2_h(k_t)$. The mappings $F^1$ and $F^2$ are stated in Eq. (62) (in \cite{ps19a}) and (\ref{discrEq:rkSchemeVecDAE}) respectively. Furthermore, matrices $F^1_{X+}$, $F^1_{X}$ are given by Eqns (70)--(71) in \cite{ps19a} and matrices $F^2_{X+}$, $F^2_{X}$ by (\ref{discrEq:FXPlusDAE})--(\ref{discrEq:FXDAE}). It remains to provide formula for $F^1_h(k_t)$ and $F^2_h(k_t)$. According to Eq. (62) (in \cite{ps19a}) and (\ref{discrEq:rkSchemeVecDAE}) we have
\begin{eqnarray}
F^1_{h}(k) = \frac{\partial F^1(X^1(k+1),X^1(k),u(k),h(k))}{\partial h(k)} & = &  \left(
\begin{array}{c}
-\sum_{j=1}^s a_{1j} f^1(x_j,u) \\  \vdots \\ 
-\sum_{j=1}^s a_{sj} f^1(x_j,u) \\  -\sum_{i=1}^s b_{i} f^1(x_i,u) 
\end{array}
\right) \nonumber \\
F^2_{h}(k) = \frac{\partial F^2(X^2(k+1),X^2(k),u(k),h(k))}{\partial h(k)} & = &  
	\left(
\begin{array}{c}
-\sum_{j=1}^s a_{1j} f^2(x_j,z_j,u) \\  0 \\\vdots \\ 
-\sum_{j=1}^s a_{sj} f^2(x_j,z_j,u) \\ 0\\  -\sum_{i=1}^s b_{i} f^2(x_i,z_i,u) 
\end{array}
\right).\nonumber
\end{eqnarray}

In order to derive analytical formula for $\pi(k_t)$ we have to figure out two vectors:  $w(k_t)=\left [F^2_{X^+}(k_t)\right ]^{-1}F^2_h(k_t)$, $w(k_t-1)=\left [F^1_{X^+}(k_t-1)\right ]^{-1}F^1_h(k_t-1)$ --- both these vectors are solutions to linear equations:

\begin{eqnarray}
&{\displaystyle 
	\left(
	\begin{array}{ccccccc}
	I - ha_{11}f_{x1}^2 & -ha_{11}f_{y1}^2 & \ldots & - ha_{1s}f_{xs}^2 & -ha_{1s}f_{ys}^2 & 0  \\
	-g_{x1}^2 & 0 & \ldots & 0 & 0 & 0   \\
	\vdots & \vdots & & \vdots & \vdots & \vdots  \\
	- ha_{s1}f_{x1}^2 & -ha_{s1}f_{y1}^2 & \ldots & I - ha_{ss}f_{xs}^2 & -ha_{ss}f_{ys}^2 & 0   \\
	0 & 0 & \ldots & -g_{xs}^2 & 0 & 0  \\
	-h b_1 f_{x1}^2 & -h b_1 f_{y1}^2 & \ldots & -h b_s f_{xs}^2 & -h b_s f_{ys}^2 & I 
	\end{array}
	\right )
	\left( 
	\begin{array}{c}
	w_{f}^1 \\ w_{g}^1 \\ \vdots \\ 
	w_{f}^s \\ w_{g}^s \\ w_f 
	\end{array}
	\right) = }\nonumber \\
&{\displaystyle 
	\left(
	\begin{array}{c}
	-\sum_{j=1}^s a_{1j} f^2(x_j,z_j,u) \\  0 \\\vdots \\ 
	-\sum_{j=1}^s a_{sj} f^2(x_j,z_j,u) \\ 0\\  -\sum_{i=1}^s b_{i} f^2(x_i,z_i,u) 
	\end{array}
	\right), }
\nonumber
\end{eqnarray}
\begin{equation}
\left(
\begin{array}{cccc}
I - ha_{11}f_{x1}^1 & \ldots & -ha_{1s}f_{xs}^1 & 0  \\
\vdots &  & \vdots & \vdots  \\
-ha_{s1}f_{x1}^1 & \ldots & I - ha_{ss}^1f_{xs} & 0  \\
-h b_1 f_{x1}^1 & \ldots & -h b_s f_{xs}^1  & I 
\end{array}
\right )
\left( 
\begin{array}{c}
w^{1} \\ \vdots \\ 
w^{s} \\  w
\end{array}
\right) = 
\left(
\begin{array}{c}
-\sum_{j=1}^s a_{1j} f^1(x_j,u) \\  \vdots \\ 
-\sum_{j=1}^s a_{sj} f^1(x_j,u) \\  -\sum_{i=1}^s b_{i} f^1(x_i,u) 
\end{array}
\right). 
\nonumber
\end{equation}

Eventually, to transform equation (\ref{adjointEquationsTrans2}) to its useful form we take into account (\ref{discrEq:discrAdjEqAdjIntDiffCalcSysEq001})--(\ref{discrEq:discrAdjEqAdjIntAlgCalcSysEq001})  and (\ref{discreteAdjSubVector}),
then the equation (\ref{adjointEquationsTrans2}) becomes
\begin{eqnarray}
&{\displaystyle \lambda^{1-}(k_t)^T \left ( h(k_t-1)\sum_{i=1}^s b_i f^1_{xi}(k_t-1) w^i(k_t-1) \right.} \nonumber\\
& {\displaystyle \left. - \sum_{i=1}^s b_i f^1(x_i(k_t),u(k_t-1))\right )   }\nonumber \\
&{\displaystyle =\lambda_f^2(k_t+1)^T \left ( h(k_t)\sum_{i=1}^s b_i f^2_{xi}(k_t) w^i_f(k_t) + h(k_t)\sum_{i=1}^s b_i f^2_{yi}(k_t) w^i_g(k_t) \right . }\nonumber \\
&{\displaystyle \left . - \sum_{i=1}^s b_i f^2(x_i(k_t+1),z_i(k_t+1),u(k_t))\right )    }\label{JumpEq3}
\end{eqnarray} 

By plugging $\lambda^{1-}(t_k)$ from equations (\ref{JumpEq3}) into equation (\ref{adjointEquationsTrans1a}) one will get
\begin{eqnarray}
&{\displaystyle \lambda_f^{2+}(k_t)^T \left ( h(k_t-1)\sum_{i=1}^s b_i f^1_{xi}(k_t-1) w^i(k_t-1) \right. }\nonumber\\ 
&{\displaystyle \left.	- \sum_{i=1}^s b_i f^1(x_i(k_t),u(k_t-1))\right )}\nonumber \\
&{\displaystyle 
-	\pi g_x(k_t) \left ( h(k_t-1)\sum_{i=1}^s b_i f^1_{xi}(k_t-1) w^i(k_t-1) \right . }\nonumber \\
&{\displaystyle \left . - \sum_{i=1}^s b_i f^1(x_i(k_t),u(k_t-1))\right )   }\nonumber \\
&{\displaystyle =\lambda^2_f(k_t+1)^T \left ( h(k_t)\sum_{i=1}^s b_i f^2_{xi}(k_t) w^i_f(k_t) + h(k_t) \sum_{i=1}^s b_i f^2_{yi}(k_t) w^i_g(k_t)- \right .} \nonumber \\
&{\displaystyle \left. \sum_{i=1}^s b_i f^2(x_i(k_t+1),z_i(k_t),u(k_t))\right )    }
\label{JumpEq4}
\end{eqnarray}
and eventually
\begin{eqnarray}
&{\displaystyle \pi(k_t) = \left [\lambda^2_f(k_t+1)^T \left (- h(k_t)\sum_{i=1}^s b_i f^2_{xi}(k_t) w^i_f(k_t) - h(k_t)\sum_{i=1}^s b_i f^2_{yi}(k_t) w^i_g(k_t)  \right . \right . }\nonumber \\
&{\displaystyle \left . +\sum_{i=1}^s b_i f^2(x_i(k_t+1),z_i(k_t),u(k_t))\right ) - \lambda_f^{2+}(k_t)^T\times}\nonumber \\
&{\displaystyle    \left . \left ( -h(k_t-1)\sum_{i=1}^s b_i f^1_{xi}(k_t-1) w_i(k_t-1) + \sum_{i=1}^s b_i f^1(x_i(k_t),u(k_t-1))\right )\right ]\bigg/   }\nonumber \\
&{\displaystyle  g_x(k_t) \left ( h(k_t-1)\sum_{i=1}^s b_i f^1_{xi}(k_t-1) w_i(k_t-1) - \sum_{i=1}^s b_i f^1(x_i(k_t),u(k_t-1))\right ). } \nonumber \\
\label{JumpEq5}
\end{eqnarray}

We assume that all functions values and their partial derivatives: $f^1$, $f^2$, $f^1_{xi}$, $f^2_{xi}$, $f^2_{yi}$ are uniformly bounded in a neighborhood of the considered solution (cf. assumptions {\bf (H1)} and {\bf (H2)} in \cite{ps19}). Furthermore, in the neighborhood, matrices $F^1_{X+}$, $F^2_{X+}$ are invertible and their elements uniformly bounded. Therefore, there exists $0<L<+\infty$ such that $\|z\|\leq L$ and $\|w\|\leq L$.

We consider RADAU IIA scheme which has the property $\sum_{i=1}^s b_i=1$. That together with the boundedness of $z$ and $w$ show that 
\begin{eqnarray}
&{\displaystyle \pi (t_k) \rightarrow \pi_t }\label{JumpEq6}
\end{eqnarray}
provided that $h(k)\rightarrow 0$.

However that property can be not sufficient to guarantee high accuracy of solutions to a continuous optimal control problem with hybrid system when a Runge--Kutta method with variable step sizes is applied. For that method we prefer a scheme for adjoint variables jumps which is determined by the formula $\pi(k_t)$ with the property 
\begin{eqnarray}
&{\displaystyle | \pi(k_t) - \pi_t| \leq O(h^{\bar{q}_{\pi}}),} \label{discEvent8}
\end{eqnarray}
where the order of convergence $\bar{q}$ can be estimated on the basis of the integration orders for state and adjoint variables.

That can be achieved for many new formula for $\pi(k_t)$ which does not dispatch far away from the formula (\ref{JumpEq5}). For example, consider the formula
\begin{eqnarray}
&{\displaystyle \hat{\pi}(k_t) = \frac{\lambda_f^{2+}(k_t)^T \left ( f^2(x(k_t),z(k_t),u(k_t-1)) - f^1(x(k_t),u(k_t-1))\right )}{g_x(k_t) f^1(x(k_t),u(k_t-1))}. }\nonumber \\
\label{discEvent9}
\end{eqnarray}

Then, we will have
\begin{eqnarray}
&{\displaystyle \hat{\pi}(t_k) - \pi_t =}\nonumber \\
&{\displaystyle g_x(x(t_t)) f^1(x(t_t^-),u(t_t^-))\lambda_f^{2+}(k_t)^T \left ( f^2(x(k_t),z(k_t),u(k_t-1)) -\right . }\nonumber \\
&{\displaystyle \left.  f^1(x(k_t),u(k_t-1))\right) \bigg/ g_x(k_t) f^1(x(k_t),u(k_t-1)){g_x(x(t_t)) f^1(x(t_t^-),u(t_t^-))}-}\nonumber\\
&{\displaystyle g_x(k_t) f^1(x(k_t),u(k_t-1))\lambda(t_t^+)^T\left (f^2(x(t_t^+),z(t_t^+),u(t_t^+))-\right . }\nonumber \\
&{\displaystyle \left . f^1(x(t_t^-),u(t_t^-))\right )\bigg/g_x(k_t) f^1(x(k_t),u(k_t-1))g_x(x(t_t)) f^1(x(t_t^-),u(t_t^-)).}\nonumber
\end{eqnarray}
Since we have to assume that in the neighborhood of $x(t_t)$ the following regularity assumption must hold (cf. the assumption {\bf (H3)} in \cite{ps19}) 
\begin{eqnarray}
&{\displaystyle \left | g_x(x) f^1(x,u(t_t))\right | \geq M > 0,}\label{discEvent11}
\end{eqnarray}
and the product of Lipschitz functions is a Lipschitz function on the bounded domain, after some transformations we will arrive at the relation
\begin{eqnarray}
\left | \hat{\pi}(t_k) - \pi_t \right | & \leq & L \left ( \|g_x(x(k_t)) - g_x(x(t_t))\| + \| \lambda_f^{2+}(k_t) - \lambda(t_t^+)\| + \right . \nonumber \\
& & \|f^1(x(k_t),u(k_t)) - f^1 (x(t_t),u(t_t)) \| + \nonumber \\
& & \left . \|f^2(x(k_t),z(k_t),u(k_t)) - f^2 (x(t_t),z(t_t),u(t_t)) \|\right ) \nonumber
\end{eqnarray}
Since 
\begin{eqnarray}
\lambda_f^{2+}(k_t) - \lambda(t_t^+) & = & O(h^{\bar{q}_f}) \nonumber \\
x(k_t) - x(t_t) & = & O(h^p) \nonumber \\
z(k_t) - z(t_t) & = & O(h^q) \nonumber
\end{eqnarray}
($ \bar{q}_f $, $ p $ and $ q $ are the appropriate integration orders) and if we assume that functions $g_x$, $f^1$, $f^2$ are Lipschitz continuous then (\ref{discEvent8}) holds with $ \bar{q}_{\pi} = min\{\bar{q}_f,p,q \} $. From (\ref{adjointEquationsTrans1a})  the adjoint variable $ \lambda^{1-}(k_t) $ is calculated with the error $ O(h^{\bar{q}_f}) $. $ \lambda^{1-}(k_t) $ states an initial condition for the calculation of the adjoint variables in the non-sliding phase. The order of the global error of the adjoint variables in the non-sliding phase is therefore reduced relative to the estimates presented in Table \ref{Orders}. The order of convergence for reduced gradients is therefore equal to $ \tilde{p}_d =2 $ for both sliding and non sliding phase.

\section{Numerical results}
\label{secResults}

In this section we present results of solving three optimal control problems with hybrid systems by using the methods discussed in this paper and in the papers accompanying it (\cite{ps19},\cite{ps19a}). The results have been obtained by our preliminary software based on its two core subroutines: the first one aimed at solving optimal control problems with piecewise constant approximations to control functions (\cite{py99}); the second one for evaluating solutions to differential--algebraic equations and their corresponding adjoint equations.

To solve the reported problems we applied the SQP code described in \cite{py99}, however instead of using range--space active set method for solving QP subproblems the new version
of the SQP code is based on the implementation of the interior point method as described in \cite{gw03}. The matrices $H_k$ used in the direction finding subproblems ${\bf P_c(u)}$ (\cite{ps19a}) were evaluated according to BFGS updates with the Powell's modifications (\cite{p77})---see the discussion in Section 4 of Chapter 5 of \cite{py99}.

We have used RADAU5 subroutine, which is the Fortran implementation of the RADAU IIA scheme. The subroutine does not have the facility for locating switching points and then restarting the procedure with a new description of differential--algebraic equations. In order to enhance RADAU5 procedure applicability to DAEs with hybrid description we incorporated into it the subroutine ROOTS (for finding a root of nonlinear algebraic equations with a secant method) from the SUNDIALS package. Since the ROOTS subroutine was implemented with multistep integration methods in mind (in which polynomial approximations of state variables $(x,y)$ between mesh points are provided by methods themselves) we had to work out an approach to the interpolation of state variables between mesh points which would be suitable for Runge--Kutta methods.

In our implementation switching points $t_t$ are evaluated with the help of  $x_i(k)$, $i=1,\ldots,s$ determined at the intermediate points $t(k)$, $t(k)+c_i h(k)$, $i=1,\ldots,s$. 
Having $c_i$ and $x_i(k)$
we evaluate derivatives at these points:
\begin{eqnarray}
&{\displaystyle      x_i^{'} (k) = f^1(x_i(k),u(k-1),t(k)+c_ih(k)),\ i=1,\ldots,s}\nonumber
\end{eqnarray}
and then construct the Hermite polynomials by taking into accounts points $x_i(k)$ and their derivatives $x_i^{'}(k)$, $i=1,\ldots,s$. In our current implementation we use two vectors (and their corresponding derivatives) to build the polynomial approximation to vector $x$ on the time interval $[t(k),t(k) + h(k)]$ (\cite{atkinson74}):
\begin{eqnarray}
x(t) & = & (2\tau^3-3\tau^2+1)x(k) + h(k) (\tau^3-2\tau^2+\tau) x^{'}(k) + \nonumber \\
& & (-2\tau^3+3\tau^2) x_s(k) + h(k) (\tau^3-\tau^2)x_s^{'}(k),\label{Sw3}
\end{eqnarray}
That interpolating scheme guarantees the accuracy $O(h(k)^4)$  (\cite{atkinson74}). We notice, that the order of the interpolating scheme is lower than the numerical integration order $ p=5 $ and affects the integration order after the discrete transition. Nevertheless, the estimated orders for reduced gradients calculation presented in Table \ref{Orders} remains valid if we assume integration order $ p=4 $.   

The use of interpolating polynomials in locating switching times can influence also the accuracy with which these times are determined. 
Suppose that on an open interval $(a,b)$, such that $t_t\in (a,b)$, the state trajectory $x(t)$ is perturbed by $\delta x(t)$. To that state perturbation corresponds the perturbation of the switching time, $t_t + \delta t_t$. We will have
\begin{eqnarray}
&{\displaystyle g(x(t_t+\delta t_t) + \delta x(t_t+\delta t_t)) = 0.}\label{discEvent12}
\end{eqnarray}
Evaluating $g$ around $(t_t,0)$ will result in
\begin{eqnarray}
0 & = & g(x(t_t+\delta t_t) + \delta x(t_t+\delta t_t)) = g_x(x(t_t))\delta x(t_t) + \nonumber \\
& & g_x(x(t_t)) f^1(x(t_t),u(t_t)) \delta t_t + o(\delta t_t, \delta x). \label{discEvent14}
\end{eqnarray}
Eventually, we have
\begin{eqnarray}
&{\displaystyle  \delta t_t = - \frac{1}{g_x(x(t_t)) f^1(x(t_t),u(t_t))}g_x(x(t_t)) \delta x(t_t) +  o(\delta t_t, \delta x).  }\label{discEvent15}
\end{eqnarray}
where $\tau= (t-t(k))/h(k)$. 

This together with (\ref{discEvent11}) imply that to the perturbation of $x$ by $O(h^p)$ corresponds the perturbation of the switching time $t_t$ of the same order $O(h^p)$. It is straightforward to show that this perturbation will have the same effect on order of convergence of integrating procedures as the perturbation of initial conditions of the same order. 

The next issue concerning the implementation of adjoint equations evaluation scheme is that related to the jump formula. To obtain the numerical results reported in the paper we have applied the formula (\ref{JumpEq5}). Since it guarantees convergence of the adjoint equations Runge--Kutta scheme to their continuous counterpart under the assumption $h(k)\rightarrow 0$ we set absolute and relative tolerances to $10^{-9}$ in RADAU5 procedure. Consequently our tolerances for satisfying optimality conditions had to be essentially less demanding. Our optimization procedure stopped iterating when the following set of conditions was satisfied: 
\begin{eqnarray}
\sigma_{c_k}^{H_k}(u_k) & \geq & -\varepsilon,  \nonumber \\
\left | g^1_i (u_k) \right | & \leq & \varepsilon,\ i\in E, \nonumber \\
g^2_j (u_k) & \leq & \varepsilon,\ j\in I.\nonumber 
\end{eqnarray}
We set $\varepsilon = 10^{-6}$ in our calculations.

We show the results of application of our numerical procedures to three optimal control problems with nonlinear differential equations. In the first two examples functions $ g(x) $ defining switching surface are linear, in the third example it is nonlinear.

{\bf Example 1} The first example concerns the Coulomb--Stribeck friction model. A mass $m$ is attached to inertial space with a spring $k$. The mass is riding on a belt, which itself is moving with a constant velocity $v_{dr}$ (see Figure 6.4a in \cite{ln04}).The relative velocity of the mass with respect to the belt is equal to $v_{rel} = v - v_{dr}$. Between the mass and the belt there is the dry friction with a friction force $F_T$. In the slip phase it is the function of $v_{rel}$ and is given by the relation
\begin{eqnarray}
&{\displaystyle F_T = -\frac{\mu_s}{1+\delta |v_{rel}|}F_N{\rm sign}(v_{rel}).}\label{FT}
\end{eqnarray}
Here, $F_N = mg$. Furthermore, in the stick phase the friction force is limited by the relation $|F_T| \leq F_s = \mu_s mg$.
The functions $f_1$, $f_2$ and $g$ of the hybrid system are 
\begin{eqnarray}
f^1(x,u) = \left [\begin{array}{c}
x_2 \nonumber \\
-\frac{k}{m}x_1 + \frac{1}{m} \frac{F_s}{1+\delta |x_2 - v_{dr}|} + x_3 \\
u \end{array} \right ], \label{f1}
\end{eqnarray} 
\begin{eqnarray}
f^2(x,u) = \left [\begin{array}{c}
x_2 \nonumber \\
-\frac{k}{m}x_1 - \frac{1}{m} \frac{F_s}{1+\delta |x_2 - v_{dr}|} +x_3 \\
u \end{array} \right ] \label{f2}
\end{eqnarray}   
and 
\begin{eqnarray}
&{\displaystyle g(x) = x_2 - v_{dr}.} \label{h}
\end{eqnarray}
Here, $x_1$ corresponds to the mass position, $x_2$ to its velocity and $x_3$ influences the mass movement through the control $u$. 

The optimal control problem is as follows
\begin{eqnarray}
&\min_{u\in {\mathcal U}} x_2(t_f)& \label{Example1}\\
&{\rm s.\ t.}\\
&\ x(t_0) = x_0& 
\end{eqnarray}
\begin{eqnarray}
x' &=& f^1(x,u),\ \text{if } g(x) < 0 \\
x' &=& f^2(x,u),\ \text{if } g(x) > 0 \\
\nonumber \\
x' &=& f_F(x,u) + g_x^T(x)z,\ \text{if sliding mode occurs} \\
0 &=& g(x) \nonumber
\end{eqnarray}
 

and 
\begin{eqnarray}
&{\displaystyle x_1(t_f) - 0.6 = 0,}
\end{eqnarray}
with the additional constraint on the control signal value $u$   
\begin{eqnarray}
&{\displaystyle  -2.5 \leq u(t) \leq 2.5\ t\in [t_0,t_f],}\nonumber
\end{eqnarray}
and the time interval endpoints $ t_0 = 0 $ $t_f=1$.

\begin{figure}
	\centering     
	\subfigure[Optimal state trajectories.]{\label{FrTrj}\includegraphics[width=0.45\textwidth]{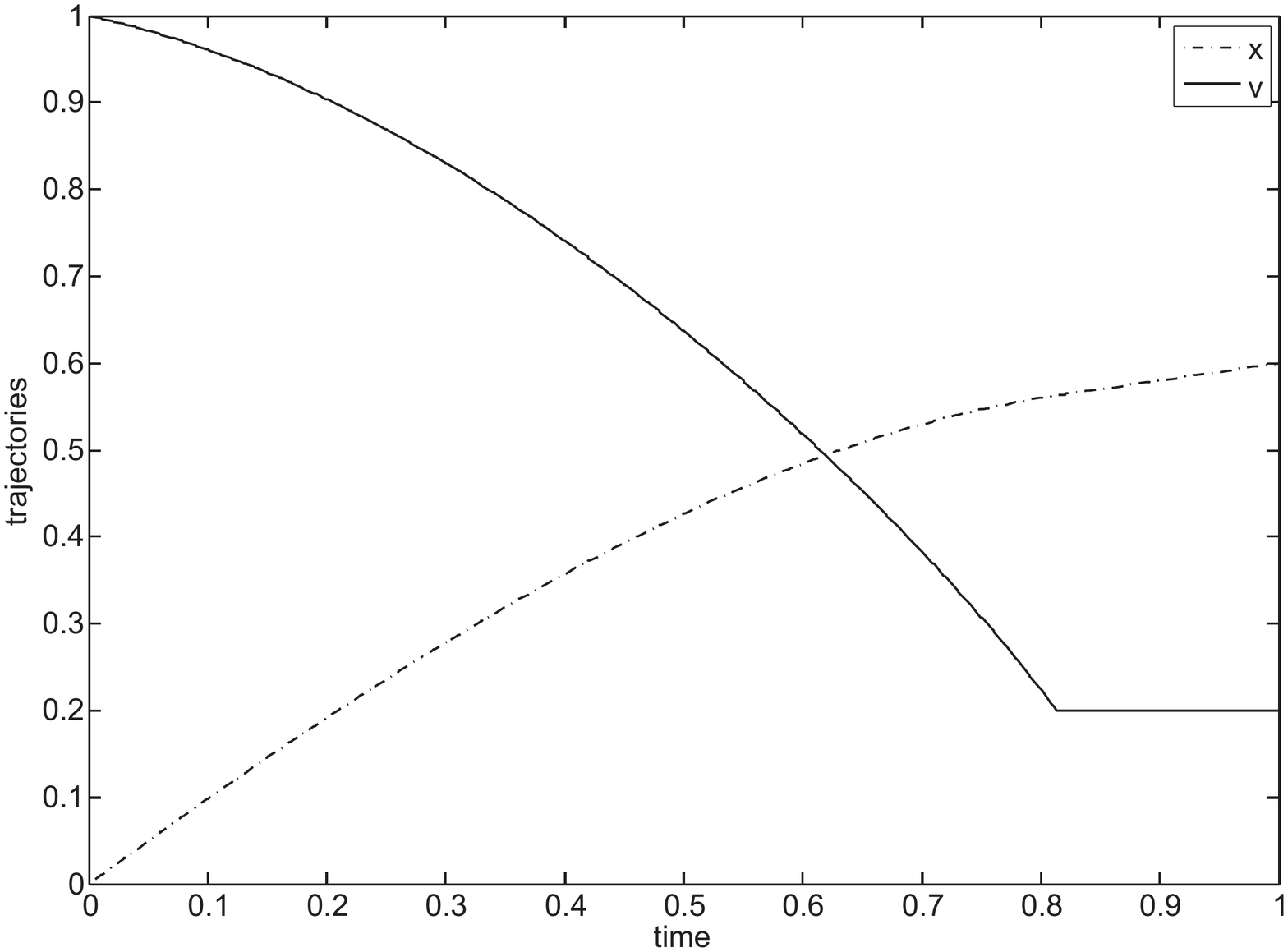}}
	\subfigure[Optimal control.]{\label{FrCtr}\includegraphics[width=0.45\textwidth]{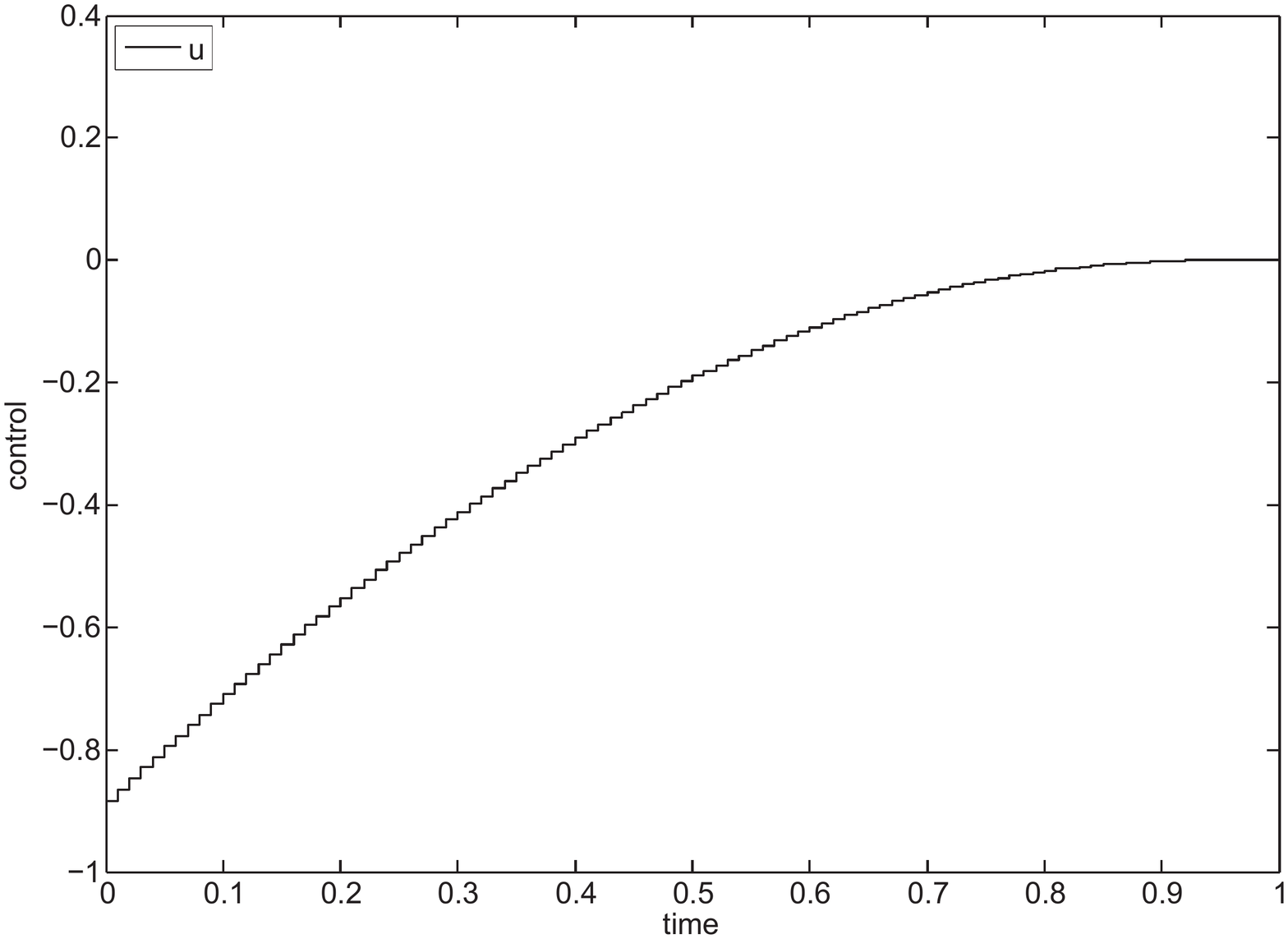}} 
	\caption{Mass--spring example optimal solution}
\end{figure}

The optimal trajectories are shown in Fig. \ref{FrTrj} and the optimal control is shown in Fig. \ref{FrCtr}. We used $N=100$ in piecewise constant approximations of control functions. One can observe that from the time $t\approx 0.8$ the system is in the sliding mode. The program needed $4$ iterations to find the solution.

{\bf Example 2}  
The example concerns the application of our approach to the problem of planning a haemodialysis process. The problem is fully described in \cite{sprn2019}. In this paper we present some results obtained for a variant of the problem. 

The system equations that determine the concentrations of urea and phosphorus in intracellular fluid -- $C_{IC}^{urea}$, $C_{IC}^{PO_4}$, urea and phosphorus concentrations in extracellular fluid -- $C_{EC}^{urea}$, $C_{EC}^{PO_4}$ and ultrafiltration volume -- $UFR$ are as follows
\begin{eqnarray}
\frac{dC_{EC}^{urea}}{dt} & = & \frac{K_{IE}^{urea} \cdot (C_{IC}^{urea}-C_{EC}^{urea}) -
	C_{EC}^{urea} \cdot (K_D^{urea}+K_r^{urea}+K^{ufr})
}{0.34 \cdot V(0)-UFR} \nonumber \\
& &  \label{dialysis1} \\
\frac{dC_{IC}^{urea}}{dt} & = & \frac{K_{IE}^{urea}\cdot(C_{EC}^{urea}-C_{IC}^{urea})+G^{urea}}{0.66 \cdot V(0)} \label{dialysis2} \\
\frac{dC_{EC}^{PO_4}}{dt} & = & \frac{K_{IE}^{PO_4} \cdot (C_{IC}^{PO_4}-C_{EC}^{PO_4}) - K_D^{PO_4} \cdot C_{EC}^{PO_4}}{0.34 \cdot V(0)-UFR} + K_{3}^{PO_4} + K_{4}^{PO_4} \label{dialysis3} \\
\frac{dC_{IC}^{PO_4}}{dt} & = & \frac{K_{IE}^{PO_4} \cdot (C_{EC}^{PO_4}-C_{IC}^{PO_4}) }{0.66 \cdot V(0)} \label{dialysis4}\\
\frac{dUFR}{dt} & = & U_{ufr} \label{dialysis5}\\
K_4^{PO_4} &=& \alpha \cdot \max \left (C_{min}^{PO_4}-C_{IC}^{PO_4},0\right ) \label{dialysis6}\\
K_3^{PO_4} &=& \beta \cdot \max \left (C_{max}^{PO_4}-C_{IC}^{PO_4},0\right ) \label{dialysis7} 
\end{eqnarray}
where
\begin{eqnarray}
&{\displaystyle K_D^{urea} = \frac{e^{K_0A \left (Q_D - Q_B\right)/\left ( Q_BQ_D\right )}-1}{e^{K_0A \left (Q_D - Q_B\right)/\left ( Q_BQ_D\right )}-Q_D/Q_B}}.\label{dialysis1a}
\end{eqnarray}
The model coefficients are explained in \cite{sprn2019}---we used the model parameters as stated for the first run of optimization. The algebraic equations are responsible for hybrid behavior of the system equations. We have $0< C_{min}^{PO_4}-C_{IC}^{PO_4}< C_{max}^{PO_4} < +\infty$ and two switching functions:
\begin{eqnarray}
g_1(C_{IC}^{PO_4}) & = & C_{IC}^{PO_4} - C_{min}^{PO_4} \label{dialysisS1} \\
g_1(C_{IC}^{PO_4}) & = & C_{IC}^{PO_4} - C_{max}^{PO_4}. \label{dialysisS2}
\end{eqnarray}

Having combined kinetic models of urea and phosphorus we look for proper concentrations of urea and phosphorus at the end of the haemodialysis process by controlling the parameters $Q_B$, $Q_D$ and $U_{ufr}$ (see \cite{sprn2019} for details). In other words, by solving the optimal control problem we want to choose a proper dialysis membrane in order to achieve final parameters of haemodialysis.

The optimization problem is as follows: 
\begin{equation}
\min_{Q_B,Q_D,U_{ufr}}C_{EC}^{urea}(t_f) \label{dialysis8}
\end{equation}
subject to the constraints (\ref{dialysis1})--(\ref{dialysis7}), the following constraints at final time 
\begin{eqnarray}
& C_{IC}^{urea}(t_f) &\leq L^{urea}_{IC} \label{dialysis10} \\
L_{min}^{UFR} \leq & UFR(t_f) &\leq L_{max}^{UFR}, \label{dialysis11}
\end{eqnarray}
and the constraints on the control variables
\begin{eqnarray}
Q_B^{min} \leq & Q_B(t) &\leq Q_B^{max} \label{dialysis12} \\
Q_D^{min} \leq & Q_D(t) &\leq Q_D^{max} \label{dialysis14} \\
U_{ufr}^{min} \leq & U_{ufr}(t) &\leq U_{ufr}^{max}  \label{dialysis15}.
\end{eqnarray}
on $[t_0,t_f]$.
Parameters for the inequality constraints are the same as in \cite{sprn2019} with the exception that now we allow $U_{ufr}$ to vary by assuming: $U_{ufr}^{min}=0$,  $U_{ufr}^{max}=75$.

\begin{figure} 
	\center{
		\includegraphics[width=10cm,height=10cm,clip,angle=-90,keepaspectratio]{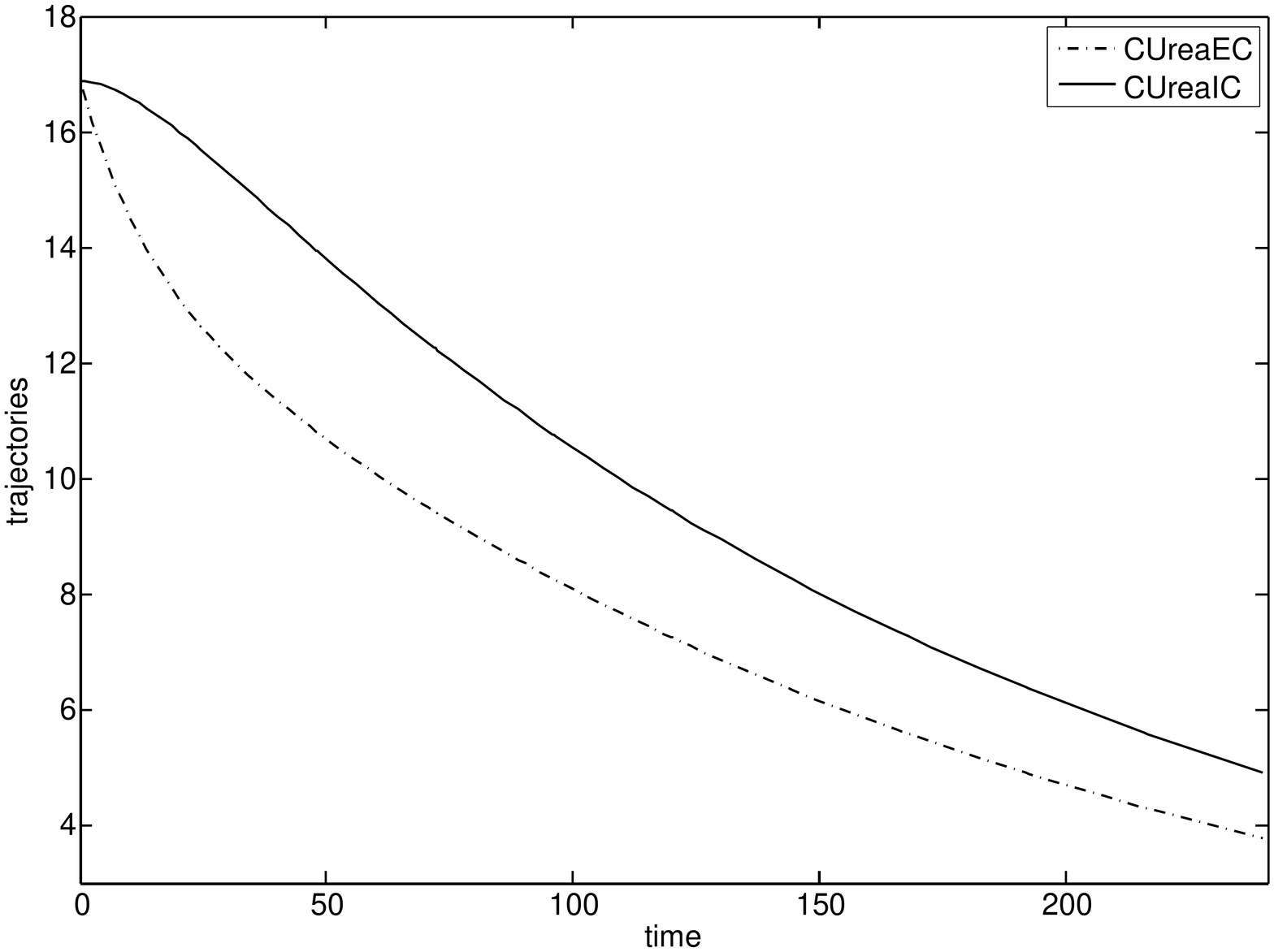}}
	\caption{State trajectories for the urea $C_{IC}^{urea}$ and $C_{EC}^{urea}$. 	} \label{DialysisTrja}
\end{figure}

\begin{figure} 
	\center{
		\includegraphics[width=10cm,height=10cm,clip,angle=-90,keepaspectratio]{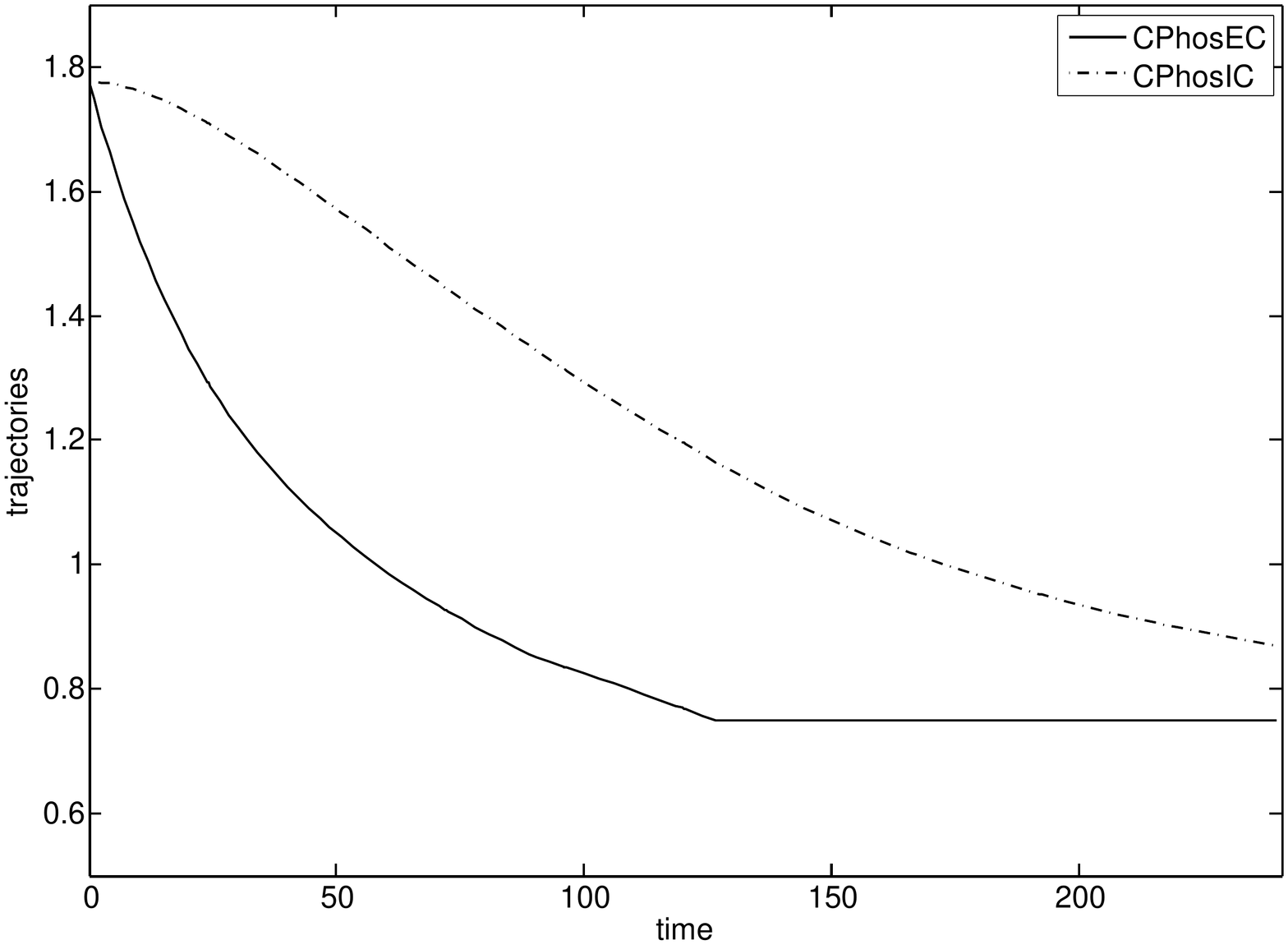}}
	\caption{State trajectories for the phosphorus $C_{IC}^{PO_4}$ and $C_{EC}^{PO_4}$ (rebound of $C_{EC}^{PO_4}$ can be seen in 125 minute of haemodialysis). 	} \label{DialysisTrjb}
\end{figure}

\begin{figure} 
	\center{
		\includegraphics[width=10cm,height=10cm,clip,angle=-90,keepaspectratio]{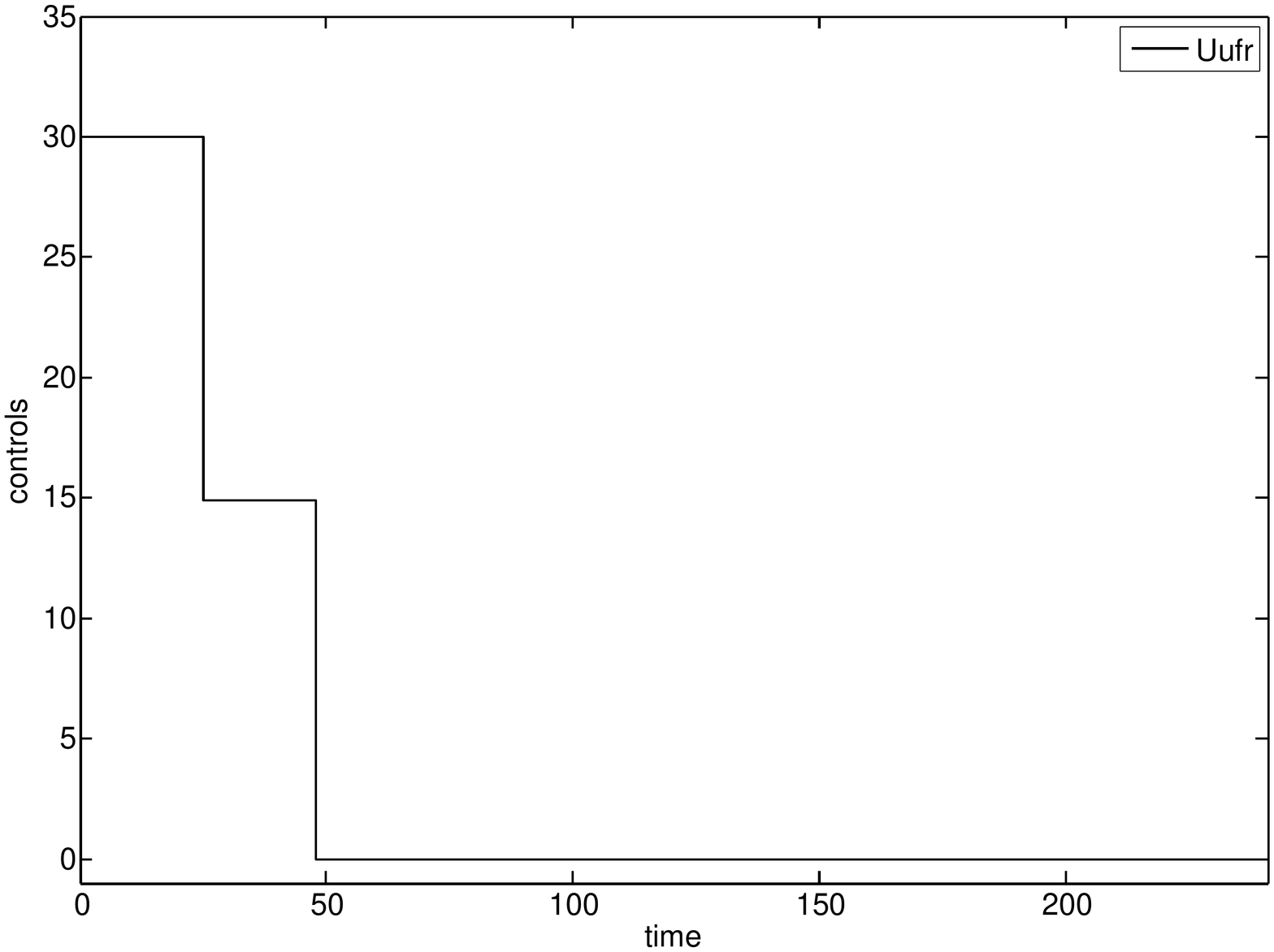}}
	\caption{Control trajectory for the $U_{ufr}$. 	} \label{DialysisCtr}
\end{figure}

Some of the optimal trajectories are shown in Figs. \ref{DialysisTrja}--\ref{DialysisTrjb}. One can observe that from the time $t\approx 125$ the system is in the sliding mode.  Optimal control variable $U_{ufr}$ is shown in Fig. \ref{DialysisCtr}, the other control variables reached their allowed maximum values---we used $N=10$ in the piecewise constant approximations of control variables. Optimization procedure needed $69$ iterations to find the solution with the specified accuracy.  

{\bf Example 3}

The example concerns the optimal control of a of a race car discussed in \cite{sa2010} and depicted in Fig.\ref{figRaceCar}.
\def\iangle{20} 
\def\down{-90}
\def\arcr{1.4cm} 
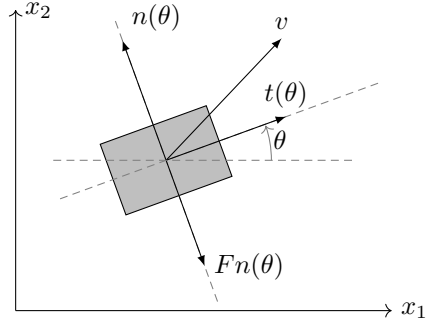
\begin{figure}[h]
	\center{
		\begin{tikzpicture}[
		force/.style={>=latex,draw=black,fill=black},
		axis/.style={densely dashed,gray,font=\small},
		M/.style={rectangle,draw,fill=lightgray,minimum size=0.5cm,thin},
		m/.style={rectangle,draw=black,fill=lightgray,minimum size=0.3cm,thin},
		plane/.style={draw=black,fill=blue!10},
		string/.style={draw=red, thick},
		pulley/.style={thick},
		]
		
		\begin{scope}[rotate=\iangle]
		\node[rectangle,draw,fill=lightgray,minimum width=1.5cm,minimum height=1.0cm,thin,transform shape] (M) {};
		
		{[axis]
			\draw (0,-2) -- (0,2) node[right] {};
			\draw (-1.5,0) -- (3,0) node[right] {};

		}
		
		{[force,->]
			\draw (0,0) -- (0,1.7) node[above right] {$n(\theta)$};
			\draw (0,0) -- (1.7,0) node[above] {$t(\theta)$};
			\draw (0,0) -- (2,1) node[above] {$v$};
			\draw (0,0) -- (0,-1.5) node[right] {$Fn(\theta)$};
		}
		
		\end{scope}
		
		\draw[axis] (-1.5,0) -- (2.5,0) node[below] {};
		\draw[->] (-2,-2) -- (-2,2) node[right] {$x_2$};
		\draw[->] (-2,-2) -- (3,-2) node[right] {$x_1$};
		\draw[axis,->,solid,shorten >=0.5pt] (0:\arcr)	arc(0:\iangle:\arcr);
		\node at (0.5*\iangle:1.1*\arcr) {$\theta$};
		
		\end{tikzpicture}
	}	
	\caption{Race car model -- mathematical description}
	\label{figRaceCar}
\end{figure}

We introduce the following notation
\begin{itemize}
	\item $ x = (x_1,x_2)^T \in \mathbb{R}^2 $ is a position of a car in the plane
	\item $ v = (v_1,v_2)^T \in \mathbb{R}^2 $ is a velocity of a car
	\item $ \theta \in \mathbb{R} $ is an angle of the car orientation
	\item $ t(\theta) = (\cos\theta,\sin\theta)^T $ is an unit vector pointing in the direction $ \theta $
	\item $ n(\theta) = (-\sin\theta,\cos\theta)^T $ is an unit vector normal to $ t(\theta) $
	\item $ a(t)\in \mathbb{R} $ is an acceleration force 
	\item $ s(t)\in \mathbb{R} $ is a steering control
\end{itemize}

The equations of motion are \cite{sa2010}
\begin{eqnarray}
x' & = & v, \label{eqRaceCarPosDer}\\
v' & = & a(t)t(\theta)+Fn(\theta), \label{eqRaceCarVelDer}\\
\theta' & = & s(t) (t(\theta)^Tv), \label{eqRaceCarThetaDer}
\end{eqnarray}
where $ Fn(\theta) $ is a friction force vector. $ F $ is given by
\begin{equation}
F = -\mu N {\rm sign}(n(\theta)^Tv), \label{eqRaceCarFricForce}
\end{equation}
where $ \mu $ and $ N $ are a friction coefficient and a normal contact force respectively. The switching surface is therefore defined as a solution set for an equation
\begin{equation}
n(\theta)^Tv = -v_1 \sin\theta  + v_2 \cos\theta = 0.
\end{equation}

During the non-sliding motion the amplitude of a friction force vector is constant and equal to $ \mu N $. The friction force is normal to $ t(\theta) $ and it acts such that it decreases the component of a velocity which is normal to a direction $ t(\theta) $. In other words the friction force always attempts to reduce the angle between  $ t(\theta) $ and $ v $ to zero. Physically, the non-sliding motion corresponds to a drift of a car.

We face the the sliding motion when the following conditions are satisfied
\begin{eqnarray}
n(\theta)^Tv & = & 0, \\
-\mu N -s(t) (t(\theta)^Tv)^2 & < & 0, \\
\mu N -s(t) (t(\theta)^Tv)^2 & > & 0. 
\end{eqnarray}
During the sliding motion the car is always oriented in the direction of a car velocity ($ t(\theta) $ and $ v $ are collinear), there is no drift and the friction force prevents the motion in a direction normal to $ t(\theta) $.

We solved the following optimization problem:
\begin{eqnarray}
&{\displaystyle \min_{a,s}x_1(t_f)}\label{Drift1}
\end{eqnarray}
subject to the system equations (\ref{eqRaceCarPosDer})--(\ref{eqRaceCarThetaDer}), the endpoint constraints
\begin{eqnarray}
x_2(t_f) & = & 0 \label{Drift2} \\
-v_1(t_f) & \leq & 0 \label{Drift3}\\
v_2(t_f) & = & 0 \label{Drift4} \\
\theta(t_f) & = & 0, \label{Drift5}
\end{eqnarray}
and the constraints on controls
\begin{eqnarray}
-0.3  \leq & a(t) &\leq 0.3 \label{Drift6}\\
-1.0 \leq & s(t) & \leq 1.0 \label{Drift7}
\end{eqnarray}
on an interval $[t_0,t_f]$ ($t_f=3.0$). We assumed $\mu N=0.5$ and initial values: $x_1(0)=0$, $x_2(0)=1$, $v_1(0)=1$, $v_2(0) = 0$.

\begin{figure}[h]
	\centering
	\includegraphics[width=10cm,height=8cm,clip]{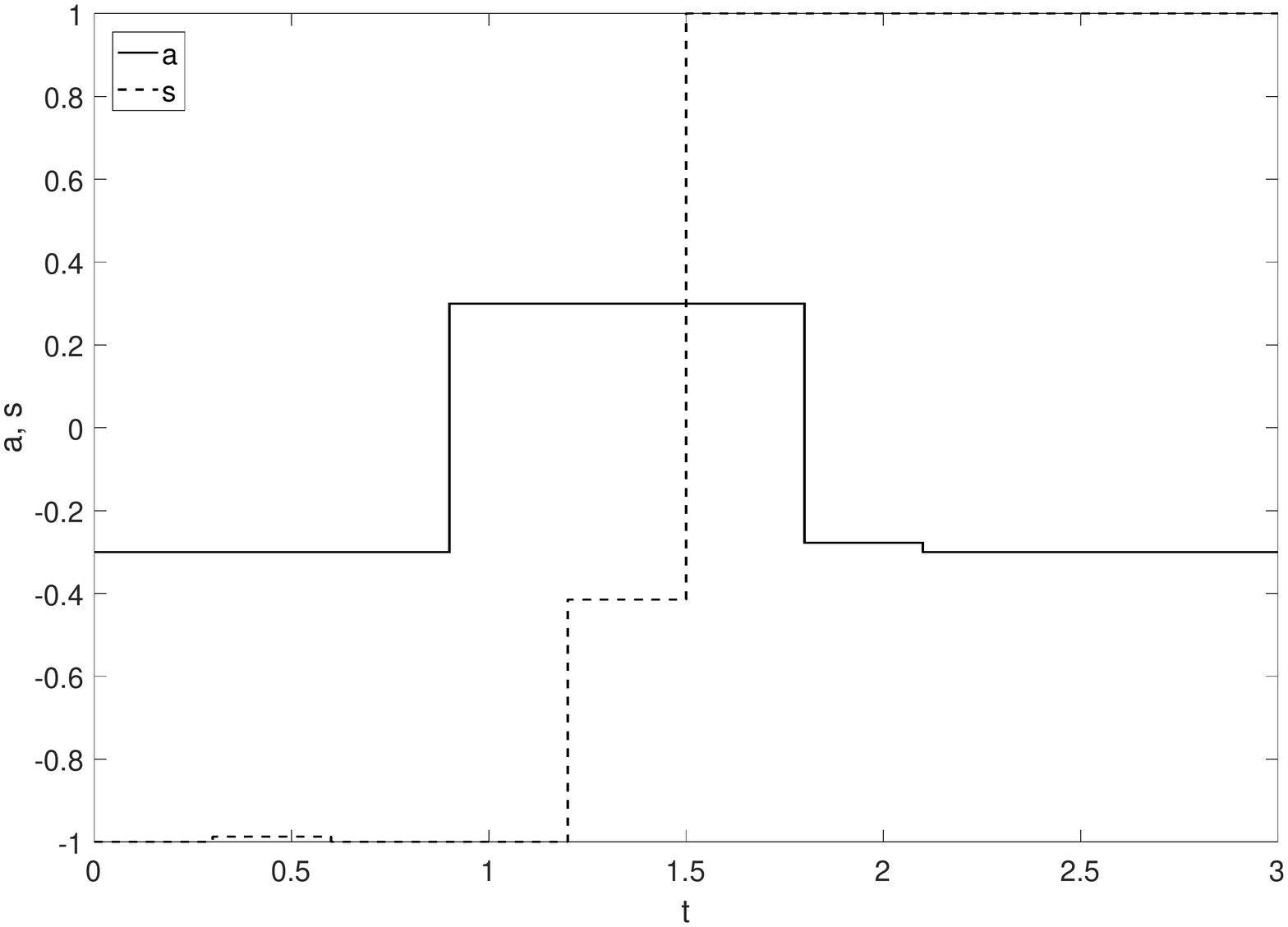}
	\caption{Optimal controls for the racing car problem.} 
	\label{figNumResRacingCarControls52Iter}
\end{figure}

\begin{figure}
	\centering     
	\subfigure[Optimal trajectories of positions.]{\label{figNumResRacingCarPositions52Iter}\includegraphics[width=0.45\textwidth]{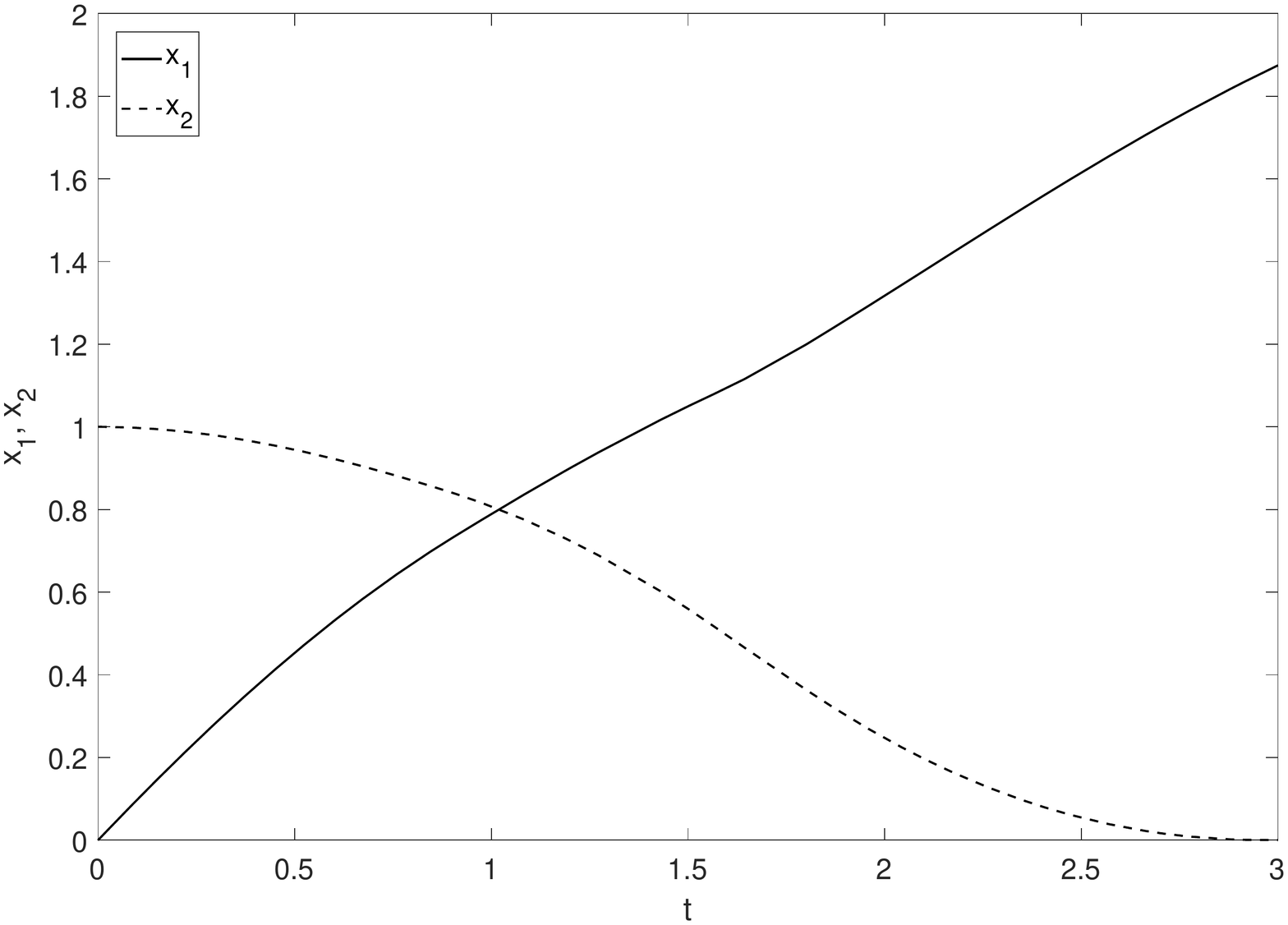}}
	\subfigure[Optimal trajectories of velocities.]{\label{figNumResRacingCarVelocities52Iter}\includegraphics[width=0.45\textwidth]{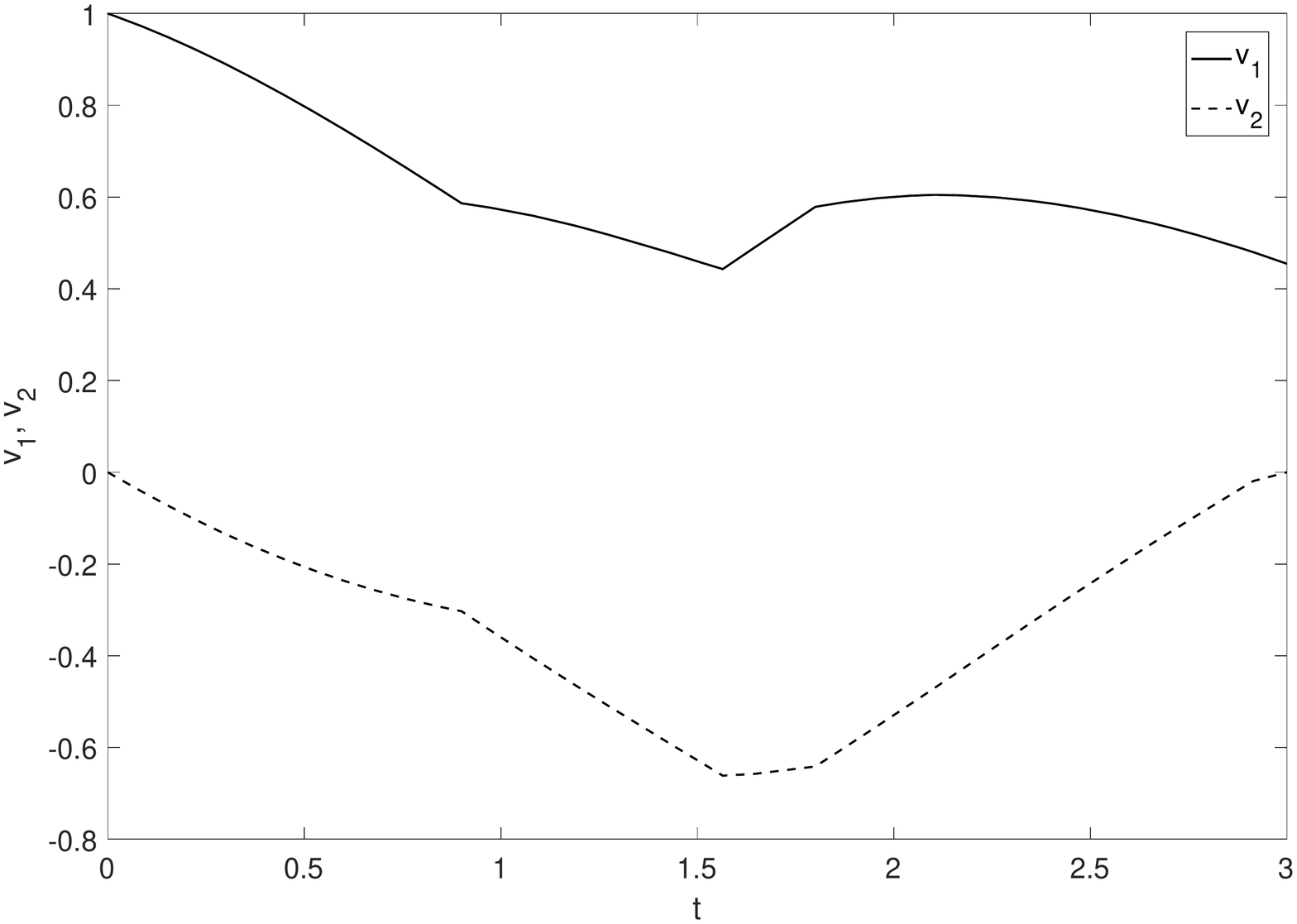}} 
	\caption{Racing car -- optimal solution}
\end{figure}

Our program needed 52 iterations to find an approximation to an optimal solution with the required accuracy, we assumed $N=10$ in the piecewise constant approximations of control functions. The optimal controls, positions and velocities are presented in Fig.\ref{figNumResRacingCarControls52Iter}, \ref{figNumResRacingCarPositions52Iter} and \ref{figNumResRacingCarVelocities52Iter} respectively.

\begin{figure}
	\centering     
	\subfigure[ $\theta $ and $ n^T(\theta)v $ trajectories.]{\label{figNumResRacingCarOrientation52Iter}\includegraphics[width=0.45\textwidth]{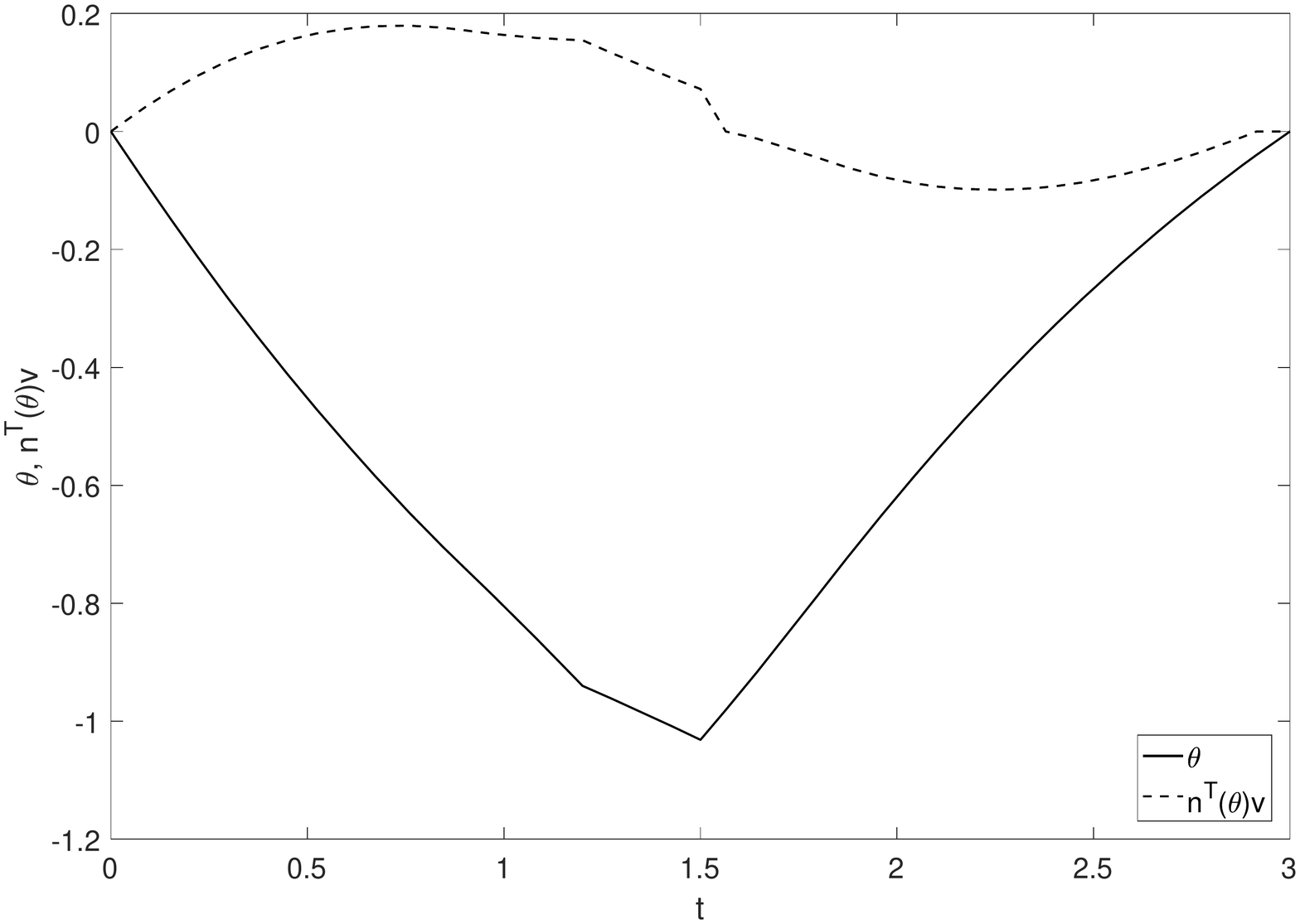}}
	\subfigure[$ \theta $ and $ n^T(\theta)v $ near the transition to the sliding mode.]{\label{figNumResRacingCarOrientationZoom52Iter}\includegraphics[width=0.45\textwidth]{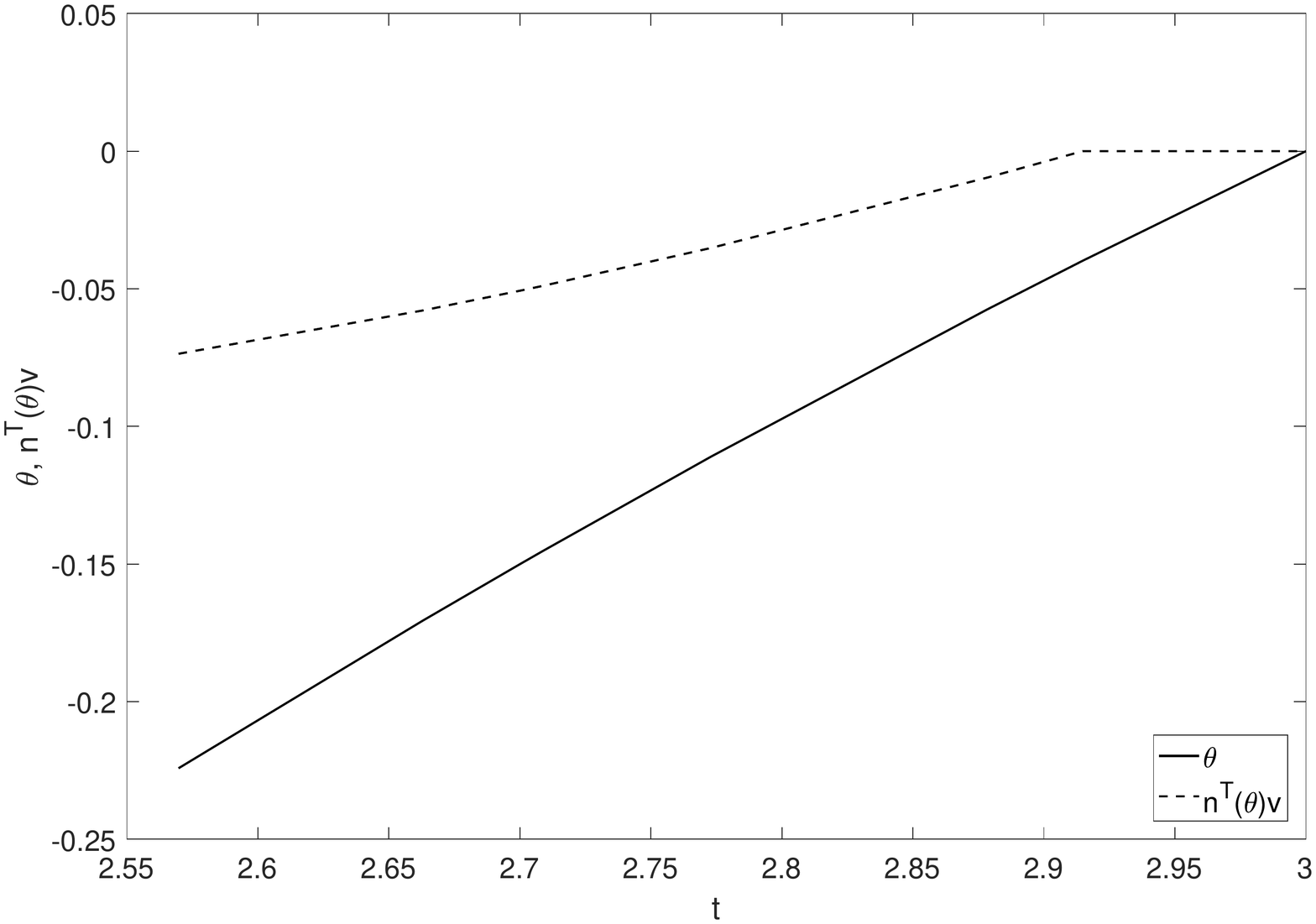}} 

	\caption{Racing car -- optimal solution}
\end{figure}

The optimal orientation trajectory as well as the value of $ n^T(\theta)v $ along the trajectory are presented in Fig.\ref{figNumResRacingCarOrientation52Iter}. 

Near the end of the time interval the condition $ n^T(\theta)v = 0$ becomes satisfied, so the system enters the sliding mode. The relevant fragment of system trajectories is presented in Fig.\ref{figNumResRacingCarOrientationZoom52Iter}.

\section{Conclusions}

The paper presents the computational approach to hybrid optimal control problems with sliding modes. It seems to be the first method for optimal control problems with hybrid systems which can exhibit sliding modes. Our computational method is based on a Runge--Kutta method for integrating system and adjoint equations. The possible improvements of our code are related to the approximating rules for the jumps of adjoint variables (such as (\ref{discEvent9})) and polynomial approximations of state variables used in the procedure of the switching times  $t_t$ localization (cf. (\ref{Sw3})). These issues are the subject of our current research.

\end{document}